\documentclass[10pt,reqno]{amsart}
\usepackage{amsmath}
\usepackage{amsfonts}
\usepackage{mathrsfs}
\usepackage{amssymb}
\usepackage{amsthm}
\usepackage{enumitem}
\usepackage{bbm}
\usepackage{times}
\usepackage{tabularx}
\usepackage{amsbsy}
\usepackage{mathtools}
\usepackage{tikz}
\usetikzlibrary{arrows,decorations.markings}
\usetikzlibrary{matrix}
\usetikzlibrary{graphs}
\usetikzlibrary{backgrounds}
\usepackage{setspace}     \spacing{1}
\usepackage{color}

\theoremstyle{Theorem}
\newtheorem{theorem}{Theorem} [section]

\newtheorem{proposition}[theorem]{Proposition}
\newtheorem{claim}[theorem]{Claim}
\newtheorem{lemma}[theorem]{Lemma}
\newtheorem{corollary}[theorem]{Corollary}

\theoremstyle{definition}
\newtheorem{definition}[theorem]{Definition}

\newtheorem{example}[theorem]{Example}
\newtheorem{remark}[theorem]{Remark}

\newtheorem*{HistoricalRemarks}{Historical Remarks}
\newtheorem{conjecture}[theorem]{Conjecture}
\theoremstyle{remark}

\newlist{enumlemma}{enumerate}{3}
\setlist[enumlemma]{label*={(\alph*)}, ref= {(\alph*)} }

\usepackage{marginnote}

\newcommand\Ts{\rule{0pt}{2.6ex}}       
\newcommand\Bs{\rule[-1.2ex]{0pt}{0pt}} 
\newcommand\Ms{\rule[-2.5ex]{0pt}{0pt}} 

\def\DynkinNodeSize{1.8mm}
\def\DynkinArrowLength{1.7mm}
\tikzset{
  dnode/.style={
    circle,
    inner sep=0pt,
    minimum size=\DynkinNodeSize,
    fill=white,
    draw},
  middlearrow/.style={
    decoration={markings,
      mark=at position 0.6 with
      {\draw (0:0mm) -- +(+135:\DynkinArrowLength); \draw (0:0mm) -- +(-135:\DynkinArrowLength); \draw (0:0mm) -- +(+315:.2pt); },
    },
    postaction={decorate}
  },
  leftrightarrow/.style={
    decoration={markings,
      mark=at position 0.999 with
      {
      \draw (0:0mm) -- +(+135:\DynkinArrowLength); \draw (0:0mm) -- +(-135:\DynkinArrowLength);
      },
      mark=at position 0.001 with
      {
      \draw (0:0mm) -- +(+45:\DynkinArrowLength); \draw (0:0mm) -- +(-45:\DynkinArrowLength);
      },
    },
    postaction={decorate}
  },
  sedge/.style={
  },
  dedge/.style={
    middlearrow,
    double distance=0.5mm,
  },
  tedge/.style={
    middlearrow,
    double distance=1.0mm+\pgflinewidth,
    postaction={draw}, 
  },
  infedge/.style={
    leftrightarrow,
    double distance=0.5mm,
  },
}

\def\ADYNK{
\begin{tikzpicture}[scale=\scales]
\useasboundingbox (-.25,-.7) rectangle (4,.5);
    \node[dnode,label=below:$\alpha_1$] (1) at (0,0) {};
    \node[dnode,label=below:$\alpha_2$] (2) at (1,0) {};
    \node[dnode,label=below:$\alpha_{\ell-1}$] (3) at (2.5,0) {};
    \node[dnode,label=below:$\alpha_\ell$] (4) at (3.5,0) {};
    \path (1) edge[sedge] (2)
          (2) edge[sedge,dashed, dash phase=1.4pt] (3)
          (3) edge[sedge] (4);
 \end{tikzpicture}
}
\def\BDYNK{
 \begin{tikzpicture}[scale=\scales]
\useasboundingbox (-.25,-.7) rectangle (4,.5);
\path           (2.45,0) edge[dedge] (3.55,0);
    \node[dnode,label=below:$\alpha_1$] (1) at (0,0) {};
    \node[dnode,label=below:$\alpha_2$] (2) at (1,0) {};
    \node[dnode,label=below:$\alpha_{\ell-1}$] (3) at (2.5,0) {};
    \node[dnode,label=below:$\alpha_\ell$] (4) at (3.5,0) {};

    \path (1) edge[sedge] (2)
          (2) edge[sedge,dashed, dash phase=1.4pt] (3)
          ;
\end{tikzpicture}
}

\def\CDYNK{
\begin{tikzpicture}[scale=\scales]
\useasboundingbox (-.25,-.7) rectangle (4,.5);
\path        (3.55,0) edge[dedge]    (2.45,0);

    \node[dnode,label=below:$\alpha_1$] (1) at (0,0) {};
    \node[dnode,label=below:$\alpha_2$] (2) at (1,0) {};
    \node[dnode,label=below:$\alpha_{\ell-1}$] (3) at (2.5,0) {};
    \node[dnode,label=below:$\alpha_\ell$] (4) at (3.5,0) {};

    \path (1) edge[sedge] (2)
          (2) edge[sedge,dashed, dash phase=1.4pt] (3)
          ;
\end{tikzpicture}}

\def\BCDYNK{
\begin{tikzpicture}[scale=\scales]
\useasboundingbox (-.25,-.7) rectangle (4,.5);
\path       (2.45,0) edge[dedge] (3.55,0); 
    \node[dnode,label=below:$\alpha_1$] (1) at (0,0) {};
    \node[dnode,label=below:$\alpha_2$] (2) at (1,0) {};
    \node[dnode,label=below:$\alpha_{\ell-1}$] (3) at (2.5,0) {};
    \node[dnode,label=below:$\alpha_\ell$] (4) at (3.5,0) {};

    \path (1) edge[sedge] (2)
          (2) edge[sedge,dashed, dash phase=1.4pt] (3)
          ;
\end{tikzpicture}
}

\def\DDYNKup{
\begin{tikzpicture}[scale=\scales]

\useasboundingbox (-.25,-1.0) rectangle (4,1.25);

    \node[dnode,label=below:$\alpha_1$] (1) at (0,0) {};
    \node[dnode,label=below:$\alpha_2$] (2) at (1,0) {};
    \node[dnode,label=below:$\alpha_{\ell-2}$] (4) at (2.5,0) {};
    \node[dnode,label=above:$\alpha_{\ell-1}$] (5) at (3.5,0.5) {};
    \node[dnode,label=below:$\alpha_\ell$] (6) at (3.5,-0.5) {};

    \path (1) edge[sedge] (2)
          (2) edge[sedge,dashed, dash phase=1.4pt] (4)
          (4) edge[sedge] (5)
              edge[sedge] (6);
\end{tikzpicture}
}

\def\EDYNK{
 \begin{tikzpicture}[scale=\scales] 
\useasboundingbox (-.25,-.7) rectangle (4,1.5);
    \node[dnode,label=below:$\alpha_1$] (1) at (0,0) {};
    \node[dnode,label=below:$\alpha_3$] (3) at (1.5,0) {};
    \node[dnode,label=right:$\alpha_6$] (4) at (1.5,1) {};
    \node[dnode,label=below:$\alpha_4$] (5) at (2.5,0) {};
    \node[dnode,label=below:$\alpha_5$] (6) at (3.5,0) {};

    \path (1) edge[sedge,dashed, dash phase=1.4pt] (3)
          (3) edge[sedge] (4)
          (3) edge[sedge] (5)
          (5) edge[sedge] (6);
\end{tikzpicture}
}
\def\EEDYNK{
 \begin{tikzpicture}[scale=\scales] 
\useasboundingbox (-.25,-.7) rectangle (4,1.5);

    \node[dnode,label=below:$\alpha_1$] (1) at (0,0) {};
    \node[dnode,label=below:$\alpha_4$] (3) at (1.5,0) {};
    \node[dnode,label=right:$\alpha_7$] (4) at (1.5,1) {};
    \node[dnode,label=below:$\alpha_5$] (5) at (2.5,0) {};
    \node[dnode,label=below:$\alpha_6$] (6) at (3.5,0) {};

    \path (1) edge[sedge,dashed, dash phase=1.4pt] (3)
          (3) edge[sedge] (4)
          (3) edge[sedge] (5)
          (5) edge[sedge] (6);
\end{tikzpicture}}

\def\EEEDYNK{
\begin{tikzpicture}[scale=\scales] 
\useasboundingbox (-.25,-.7) rectangle (4,1.5);

    \node[dnode,label=below:$\alpha_1$] (1) at (0,0) {};
    \node[dnode,label=below:$\alpha_5$] (3) at (1.5,0) {};
    \node[dnode,label=right:$\alpha_8$] (4) at (1.5,1) {};
    \node[dnode,label=below:$\alpha_6$] (5) at (2.5,0) {};
    \node[dnode,label=below:$\alpha_7$] (6) at (3.5,0) {};

    \path (1) edge[sedge,dashed, dash phase=1.4pt] (3)
          (3) edge[sedge] (4)
          (3) edge[sedge] (5)
          (5) edge[sedge] (6);
\end{tikzpicture}
}
\def\FDYNK{
\begin{tikzpicture}[scale=\scales]
\useasboundingbox (-.25,-.7) rectangle (4,.5);
\path (.95,0) edge[dedge] (2.05,0);
    \node[dnode,label=below:$\alpha_1$] (1) at (0,0) {};
    \node[dnode,label=below:$\alpha_2$] (2) at (1,0) {};
    \node[dnode,label=below:$\alpha_3$] (3) at (2,0) {};
    \node[dnode,label=below:$\alpha_4$] (4) at (3,0) {};

    \path (1) edge[sedge] (2)
          (3) edge[sedge] (4)
          ;
\end{tikzpicture}
}
\def\GDYNK{
          \begin{tikzpicture}[scale=\scales]
\useasboundingbox (-.25,-.7) rectangle (4,.5);
  \path (-.05,0) edge[tedge] (1.05,0);
    \node[dnode,label=below:$\alpha_1$] (1) at (0,0) {};
    \node[dnode,label=below:$\alpha_2$] (2) at (1,0) {};

\end{tikzpicture}
}

\def\DynkTable{
\begin{table}[h]
\footnotesize
\def\scales{.7}
\caption{Roots systems, highest and 2nd highest roots, and resonant codimension of maximal parabolic subalgebras}
\label{table:bull}
\begin{center}
\begin{tabular}{|c|m{83pt} |l|}
\hline
& Dynkin diagram and
& Highest root  $\delta$ and second-highest root  $\delta'$; 
\\
&   simple roots 
&   resonant codimension $\bar r(\lieq_{j})$ where $\lieq_j = \lieq_{\Pi\sm \{ \alpha_j\}}$
\\    \hline

$A_\ell$&
\ADYNK
&
\begin{tabular}{l  }

$ \delta = \alpha_1 + \dots + \alpha _\ell$  \Ts  \Ms
 \\
 $\bar r(\lieq_{j})= \frac1 2 
\big(
(\ell+1)^2 - j^2 - (\ell+1-j)^2 
\big)$\Bs 
 \end{tabular}
 \\		\hline
$B_\ell$&
\BDYNK
&
 \begin{tabular}{l}
$\delta =  \alpha _1 + 2 \alpha _2 + \dots + 2 \alpha _\ell  $\Ts  \Ms
\\
$\bar r(\lieq_{j})= \frac1 2 
\big(
\ell (2\ell+1)  - j^2 - (\ell-j) (2(\ell-j)+1)
\big)$\Bs
\\
\end{tabular}
\\      \hline
$C_\ell$&
\CDYNK
&
\begin{tabular}{l l }
$\delta  = 2\alpha_1 + 2\alpha_2 + \dots + 2 \alpha_{\ell-1} + \alpha _ \ell   $\Ts   \\
$\delta'  = \alpha_1 +  2\alpha_2 + \dots + 2 \alpha_{\ell-1} + \alpha _ \ell  $\Ms
\\
$\bar r(\lieq_{j})= \frac1 2 
\big(
\ell (2\ell+1)  - j^2 - (\ell-j) (2(\ell-j)+1)
\big)$\Bs
\\
\end{tabular}
\\		\hline
$BC_\ell$&
\BCDYNK
&
\begin{tabular}{l l }
$ \delta= 2 \alpha_1 + 2\alpha _{2 } + \dots + 2\alpha _{\ell -1 }  +2 \alpha _\ell $ \Ts\\
$ \delta' = \alpha_1 + 2\alpha _{2 }  + \dots + 2\alpha _{\ell -1 } + 2\alpha _\ell   $\Bs
\end{tabular}
\\ 		\hline
$D_\ell$&
\DDYNKup
&
\begin{tabular}{  l  }
$ \delta = \alpha_1 + 2\alpha _2 +\dots + 2\alpha _{\ell -2 } +   \alpha _{\ell-1}+  \alpha _\ell  $ \Ts \Ms\\
$\bar r(\lieq_{j})= \frac1 2 
\big(
 \ell (2\ell-1)  - j^2 - (\ell-j) (2(\ell-j)-1) 
\big)$
\\
\hfill for $\phantom{\ell}1\le j\le \ell-2$  \\
$\bar r(\lieq_{j})= \frac1 2 
\big(
 \ell (2\ell-1) - \ell^2
\big)$
\hfill   \Bs
 for  $\phantom{1}\ell-1\le j\le \ell$
\end{tabular}
\\		 \hline

$E_6$&
\EDYNK
&
  \begin{tabular}{l}
$ \delta =  \alpha _1 + 2\alpha _2 + 3 \alpha _3 + 2\alpha _4 +  \alpha _5 + 2  \alpha _6$\Ts \Ms \\
$\bar r(\lieq_{1})= 16 \quad
\bar r(\lieq_{2})= 25 \quad 
\bar r(\lieq_{3})= 29 \quad  $\\$
\bar r(\lieq_{4})= 26\quad
\bar r(\lieq_{5})= 16\quad 
\bar r(\lieq_{6})= 21$\Bs
\end{tabular}
  \\ 		\hline
$E_7$&
\EEDYNK
&
 \begin{tabular}{l}
$\delta = \alpha_1+2\alpha_2+3\alpha_3+4\alpha_4+3\alpha_5+2\alpha_6+ 2\alpha_7$\Ms \Ts \\
$\bar r(\lieq_{1})= 27\quad 
\bar r(\lieq_{2})= 42 \quad
\bar r(\lieq_{3})= 50\quad  $\\$
\bar r(\lieq_{4})= 53\quad 
\bar r(\lieq_{5})= 47\quad 
\bar r(\lieq_{6})= 33\quad  $\\$
\bar r(\lieq_{7})= 42$\Bs
 \end{tabular}
\\ 			
\hline
$E_8$&
\EEEDYNK
&
 \begin{tabular}{l}
 $ \delta = 2 \alpha _1 + 3\alpha _2 + 4 \alpha _3 + 5\alpha _4 + 6\alpha _5 +   $\Ts \\
\quad \quad  \quad\quad $ 4\alpha _6 + 2\alpha _7+3\alpha_8$\\
 $ \delta' =  \alpha _1 + 3\alpha _2 + 4 \alpha _3 + 5\alpha _4 + 6\alpha _5 + $ \\
\quad \quad  \quad\quad $   4\alpha _6 + 2\alpha _7+3\alpha_8$\Ms 
\\
$\bar r(\lieq_{1})= 57\phantom{0}\quad 
\bar r(\lieq_{2})= 83\phantom{0} \quad
\bar r(\lieq_{3})= 97\phantom{0}\quad  $\\$
\bar r(\lieq_{4})= 105\quad 
\bar r(\lieq_{5})= 106\quad 
\bar r(\lieq_{6})= 98\phantom{0}\quad  $\\$
\bar r(\lieq_{7})= 78\phantom{0}\quad
\bar r(\lieq_{8})= 92$\Bs
 \end{tabular}
\\ 			\hline
$F_4$&
\FDYNK  &
\begin{tabular}{l  }
$ \delta =  2\alpha _1 + 3\alpha _2 + 4 \alpha _3 + 2\alpha _4 $\Ts \\
$ \delta' =  \alpha _1 + 3\alpha _2 + 4 \alpha _3 + 2\alpha _4 $\Ms
\\
$\bar r(\lieq_{1})= 15 \quad
\bar r(\lieq_{2})= 20 \quad 
\bar r(\lieq_{3})= 20 \quad  
\bar r(\lieq_{4})= 15$\Bs
 \end{tabular}
\\ 			\hline
$G_2$&
\GDYNK  &
\begin{tabular}{l   }
$ \delta = 2 \alpha _1 + 3\alpha _2 $\Ts\Ms \quad \quad \quad
$ \delta'  =  \alpha _1 + 3\alpha _2 $
\\
$\bar r(\lieq_{1})= 5 \phantom{0}\quad
\bar r(\lieq_{2})= 5\phantom{0}$\Bs
\end{tabular}
\\ 			\hline
\end{tabular}
\end{center}
\label{default}
\end{table}
}

\def\RP{\R P}

\newcommand{\restrict}[2]{{#1}{|_{{ #2}}}}
\newcommand{\diff}{\mathrm{Diff}}

\renewcommand{\epsilon}{\varepsilon}

\def \ad{\mathrm{ad}}
\def \diag{\mathrm{diag}}
\def \Ad{\mathrm{Ad}}

\newcommand{\wtd}{\widetilde}

\DeclareMathOperator{\Diff}{Diff}
\DeclareMathOperator{\Isom}{Isom}
\DeclareMathOperator{\Prob}{Prob}
\newcommand{\sm}{\smallsetminus}
\newcommand{\R}{\mathbb {R}}

\newcommand{\Z}{\mathbb {Z}}
\newcommand{\N}{\mathbb {N}}
\newcommand{\T}{\mathbb {T}}
\newcommand{\C}{\mathbb {C}}

\newcommand{\e}{\epsilon}

\def \sl{\mathfrak{sl}}
\newcommand{\Gl}{\mathrm{GL}}
\newcommand{\Sl}{\mathrm{SL}}
\newcommand{\PSl}{\mathrm{PSL}}
\newcommand{\Sp}{\mathrm{Sp}}
\newcommand{\So}{\mathrm{SO}}
\def\SO{\So}
\newcommand{\inv}{^{-1}}
\newcommand{\id}{\mathrm{Id}}

\def\calA{\mathcal A}
\def\calE{\mathcal E}

\def\calL{\mathcal L}

\def \RP{\R P}

\def\orb{\mathcal O}

\def\ae{a.e.\ }

 \newcommand{\lieg}{\mathfrak g}
\newcommand{\lieh}{\mathfrak h}
\newcommand{\liek}{\mathfrak k}

\newcommand{\liec}{\mathfrak c}
\newcommand{\liel}{\mathfrak l}
\newcommand{\lien}{\mathfrak n}
\newcommand{\liea}{\mathfrak a}

\newcommand{\liep}{\mathfrak p}
\newcommand{\lieu}{\mathfrak u}
 \newcommand{\lieq}{\mathfrak q}

\def\wl{\mathrm{l}}
\def\Folner{F{\o}lner }
\def\vecm{\mathbf{m}}

\def\blue{\color{\blue}}

\keywords{Zimmer program, actions of lattices, lattices in semisimple Lie groups,  actions of abelian groups, Ratner theory, property (T), Lyapunov exponents, measure rigidity}


\title[Zimmer's conjecture]{Zimmer's conjecture: Subexponential growth, measure rigidity, and  strong property (T)}

\author[A.~Brown]{Aaron Brown}
\address{University of Chicago, Chicago, IL 60637, USA}
\email{awb@uchicago.edu}

\author[D.~Fisher]{David Fisher}
\address{Indiana University, Bloomington, Bloomington, IN 47401, USA}
\email{fisherdm@indiana.edu}

\author[S.~Hurtado]{Sebastian Hurtado}
\address{University of Chicago, Chicago, IL 60637, USA}
\email{shurtados@uchicago.edu}

\thanks{D.~F. was partially supported by NSF Grant DMS-1308291.  This work was begun while
he visited Chicago, a visit partially supported by NSF RTG Grant DMS-1344997.   Travel for A.~B. was supported by NSF grants DMS 1107452, 1107263, 1107367, ``RNMS: Geometric Structures and Representation Varieties" (the GEAR Network).}

\long\def\symbolfootnote[#1]#2{\begingroup\def\thefootnote{\fnsymbol{footnote}}
\footnote[#1]{#2}\endgroup}

\begin{document}
\maketitle
\begin{abstract}
We prove several cases of Zimmer's conjecture for actions of higher-rank, cocompact lattices  on low-dimensional manifolds. For example, if $\Gamma$ is a cocompact lattice in $\Sl(n, \mathbb R)$, $M$ is a compact manifold, and $\omega$ a volume form on $M$ we show that any homomorphism $\alpha\colon \Gamma \rightarrow \Diff(M)$ has finite image if the dimension of $M$ is less than $n-1$ and that any homomorphism $\alpha\colon \Gamma \rightarrow \Diff(M,\omega)$ has finite image if the dimension of $M$ is less than $n$.
The key step in  the proof is to show that any such action has uniform subexponential growth of derivatives.
This is established using  ideas  from the smooth ergodic theory of higher-rank abelian groups, structure theory of semisimple groups, and results from homogeneous dynamics.
Having established uniform subexponential growth of derivatives, we apply Lafforgue's strong property (T) to establish the existence of an invariant Riemannian metric.
\end{abstract}

\section{Introduction}

\subsection{Results, history, and motivation.}
As a special  case of our main result, Theorem \ref{theorem:main} below,  we confirm   Zimmer's conjecture for actions of cocompact lattices in $\Sl(n,\R)$.

\begin{theorem}\label{slnr}
For $n\ge 3$, let $\Gamma< \Sl(n,\mathbb R)$ be a cocompact lattice. Let $M$ be a compact manifold. If $\dim(M)< n-1$
then any homomorphism $\Gamma\rightarrow \Diff^{2}(M)$ has finite image.  In addition, if $\omega$ is a volume form on $M$
and $\dim(M)=n-1$ then any homomorphism $\Gamma \rightarrow \Diff^{2}(M,\omega)$ has finite image.
\end{theorem}


The key step in the proof is to establish that the derivatives of group elements for such an action grow subexponentially relative to their word length.  This is inspired by the third author's paper on the Burnside problem for diffeomorphism groups \cite{Hurtado}.  To prove subexponential growth of derivatives in this context, we study the induced $G$-action on a suspension space and apply a number of measure rigidity results including Ratner's theorem and  recent work of the first author with Rodriguez Hertz and Wang.  Having   established subexponential growth of derivatives, the main theorem is established by using the  {strong Banach property (T)} of Lafforgue to find an invariant Riemannian metric. The proof has many surprising features, including  its use of hyperbolic dynamics to prove an
essentially elliptic result and   its use of  results from  homogeneous dynamics to prove results about non-linear actions. We include a detailed sketch of the proof at the end of the introduction.

Theorem \ref{slnr} lies in the context of the {Zimmer Program}. In \cite{MR682830} Zimmer made a number of conjectures concerning smooth volume-preserving  actions of lattices in higher-rank semisimple groups on low-dimensional manifolds.
 These conjectures were clarified in \cite{MR934329,MR900826} and extended to the case of  smooth non-volume-preserving actions  by Farb and Shalen in \cite{MR1666834}.


The Zimmer program is motivated by earlier results on rigidity of  linear representations of lattices in higher-rank Lie groups.  The history of the subject begins in the early 1960s with results of Selberg and Weil which established that cocompact lattices in simple Lie groups other than $\PSl(2,\R)$ were locally rigid:  any perturbation of a lattice is  given by conjugation by a small group element \cite{Selberg, Weil-I}.  In the late 60s and early 70s, this was  improved by Mostow to  a global rigidity theorem showing that any isomorphism between cocompact lattices in the same class of groups extended to an isomorphism of the ambient Lie group \cite{Mostow-book}.  The global rigidity result was extended by Margulis and Prasad to include non-uniform lattices
\cite{Margulis-nonuniformtwo, MR0385005}.  These developments led to Margulis' work on superrigidity and arithmeticity in which Margulis classified all linear representations of irreducible lattices in Lie groups of higher real rank  \cite{MR1090825} and established that all such  lattices are arithmetic.

Inspired by Margulis' superrigidity theorem, in the early 1980s Zimmer proved a superrigidity theorem for cocycles from which he   proved results about orbit equivalence of higher-rank group actions \cite{MR595205}.  Motivated by earlier results in the rigidity of linear representations and the cocycle superrigidity theorem, Zimmer proposed studying non-linear representations of lattices in higher-rank simple Lie groups.  That is, given a lattice $\Gamma\subset G$, rather than studying linear representations $\rho\colon \Gamma \to \Gl(d,\R)$, Zimmer proposed studying representations $\alpha \colon \Gamma \to\Diff(M)$ where $M$ is a compact manifold.  The main objective of the Zimmer program is to show that all such non-linear representations $\alpha$ are of an ``algebraic origin.''  In particular, the  Zimmer conjecture states that if the dimension of $M$ is sufficiently small (relative to data associated to $G$) then any action $\alpha \colon \Gamma \to\Diff(M)$ should preserve a smooth Riemannian metric and thus  factor through the action of a finite group under certain additional dimension constraints.  See Conjecture \ref{conjecture:zimmer} for a precise formulation.

In this paper we establish the non-volume-preserving case of Zimmer's conjecture for actions of cocompact lattices in higher-rank split simple Lie groups as well as certain volume-preserving cases.    While there have been a number of sharp results for  actions on extremely low-dimensional manifolds   (for manifolds of dimension $1$ or $2$) or under  strong regularity conditions on the action or algebraic conditions on the lattice, prior to this paper the exact result conjectured by Zimmer was only known for non-uniform lattices in $\Sl(3,\R)$.
Our results provide  a class of higher-rank Lie groups and a large collection of lattices such that the critical dimension is as conjectured in the non-volume-preserving and either as conjectured or almost as conjectured in the volume-preserving case. In addition to establishing the conjecture for cocompact lattices in split simple Lie groups,  we also give strong partial results for actions of cocompact lattices in non-split simple Lie groups.

In the case of  volume-preserving actions, the conjecture is motivated by the following  corollary of  Zimmer's cocycle superrigidity theorem:    all volume-preserving actions in sufficiently low dimensions  preserve a  measurable Riemannian metric \cite{MR595205}. From this point of view, the main step in proving the conjecture is to promote a measurable metric to a  smooth metric.   Conditional and partial results  verifying the existence of a smooth invariant metric in the volume-preserving case are contained in many papers of Zimmer of which \cite{MR900826} provides an excellent overview.

Perhaps the best evidence for the conjecture in the case of    volume-preserving actions is Zimmer's result that all actions satisfying the conjecture have discrete spectrum \cite{MR1147291}.
 In the  non-volume-preserving case, evidence for the conjecture follows from   the works of Ghys and of Farb and Shalen on analytic actions and  work of Nevo and Zimmer that produces {measurable projective quotients} for actions which do not preserve a measure \cite{MR1254981,MR1666834, MR1933077}. 

Other strong evidence for the conjectures is provided by a plethora of results concerning actions on compact
manifolds of dimension $1$ or $2$.  The earliest results were those of Witte Morris proving that all $C^0$ actions
on $S^1$ of $\Sl(n,\mathbb Z)$ and $\Sp(2n, \mathbb Z)$ and their finite-index subgroups factor through finite
groups \cite{MR1198459}.  Later results of Burger and Monod and of Ghys show similar results for $C^1$ actions of  all lattices in higher-rank simple Lie groups \cite{MR1911660, MR1703323}.
Ghys' result also includes results for irreducible lattices in products of rank-1 groups, which admit infinite actions on the circle.
In dimension 2,  results of  Polterovich and of Franks and Handel show that all volume-preserving actions of non-uniform lattices on surfaces are also all finite \cite{MR2219247, MR1946555}.  Moreover, Franks and Handel showed that for any surface of genus at least 1, any action by a  non-uniform  lattice in a higher-rank simple Lie group  which  preserves a Borel  probability measure is finite.    Some earlier results on actions on surfaces, such as those of Farb and Shalen in the analytic category, do not require an invariant measure but instead make stronger assumptions on the acting group and the regularity of the action.  Combined with results of \cite{MR2219247} and \cite{AWBFRHZW-latticemeasure}, we resolve the conjecture almost completely for $C^2$-actions on surfaces of genus at least 1 in Theorem \ref{thm:easycor2}.
Above dimension $2$, very little is known.
See the second author's survey of the Zimmer program  \cite{F11} for a detailed history as well as earlier surveys by Feres and Katok, Labourie, and Witte Morris and Zimmer  \cite{MR1928526,MR1648087,MR2457556}.

We recall the full conjecture of Zimmer as extended by Farb and Shalen.  Given a semisimple
Lie group $G$,   let $n(G)$ denote the minimal dimension of a non-trivial real representation of the Lie algebra $\lieg$ of  $G$ and let $v(G)$ denote the minimal codimension of a maximal (proper) parabolic subgroup $Q$ of $G$.
Let $d(G)$ denote the minimal dimension of all non-trivial homogenous spaces $K/C$ as  $K$ varies over all compact real-forms of all simple factors of  the complexification of $G$.

\begin{conjecture}[Zimmer's Conjecture]\label{conjecture:zimmer}
Let $G$ be a  connected, semisimple Lie group with finite center, all of whose almost-simple factors have real-rank at least $2$.  Let $\Gamma<G$ be
a lattice.  Let $M$ be a compact manifold and let $\omega$ be a volume form on $M$.  Then
\begin{enumerate}
\item if $\dim(M) < \min (n(G), d(G), v(G))$ then any  homomorphism $\alpha\colon \Gamma \rightarrow \Diff(M)$ has finite image;
\item if $\dim(M) < \min(n(G), d(G))$ then any homomorphism $\alpha\colon \Gamma \rightarrow \Diff(M, \omega)$ has finite image;
\item if $\dim(M) < n(G)$ then for any homomorphism $\alpha\colon \Gamma \rightarrow \Diff(M, \omega)$, the image $\alpha(\Gamma)$ preserves a Riemannian metric;
\item if $\dim(M) < v(G)$ then for any homomorphism $\alpha\colon \Gamma \rightarrow \Diff(M)$, the image $\alpha(\Gamma)$ preserves a Riemannian metric.
\end{enumerate}
\end{conjecture}

Theorem \ref{slnr} verifies the conjecture for cocompact lattices in $\Sl(n,\R)$; we will discuss other cases below.
The conjecture is almost sharp in several senses.  In dimension $v(G)$, any subgroup of $G$ admits an infinite image, non-isometric, non-volume-preserving action in dimension $v(G)$, namely, the projective left-action on $G/Q$ where $Q$ is a parabolic subgroup of codimension $v(G)$.   These actions are the natural analogue of the action of $\Sl(n,\R)$ and its lattices on $\RP^{n-1}$.  In dimension $n(G)$, there is always a semisimple Lie group with finite center $\hat G$ with the same Lie algebra as $G$, a lattice $\Gamma\subset G$, and a volume-preserving, non-isometric action   on the compact manifold $\T^{n(G)}$.  However, in these examples the lattice $\Gamma$ is, in fact, the integer points of $\hat G$ with respect to the rational structure for which the representation in dimension $n(G)$ is rational; in particular, in such examples $\Gamma$ is necessarily  non-uniform.  This construction is the natural analogue of the action of $\Sl(n, \Z)$ on $\mathbb T^n$.   In particular,  $n(G)$ is a sharp  bound for results about actions of all lattices in a Lie group $G$  but  may not be sharp for results about   actions of a particular  lattice; given our results it is natural to ask if sharper bounds can be established for cocompact lattices.  Lastly, the number $d(G)$  bounds the dimension in which infinite isometric actions can occur.  
The existence of an invariant Riemannian  metric $g$ for the action $\alpha$  implies that the action is given by a homomorphism $\alpha \colon \Gamma \rightarrow K$ where $K = \Isom(M,g)$ is a compact Lie group; see discussion in Section \ref{section:finalarguments} below.  Margulis' superrigidity theorem implies that $\alpha(\Gamma)$ cannot be infinite below dimension $d(G)$.  In fact,  in the presence of an invariant metric for low dimensional actions, Margulis' superrigidity theorem classifies the possible isometry groups and  elementary geometry gives sharper results on manifolds admitting infinite, isometric actions.

\begin{HistoricalRemarks} Items $(2)$ and $(3)$ are due to Zimmer.
Zimmer stated $(2)$ in slightly different terms that were not sharp.
Item $(1)$ is a natural extension by Farb-Shalen. The conjecture as stated in both \cite{MR1666834, F11} assumed erroneously that one always has $v(G) = n(G)-1$ so the conjecture is slightly misstated in those references.  Item $(4)$ is new here, but is a natural extension of the other conjectures. We are intentionally vague concerning regularity of the diffeomorphisms in the conjecture. Zimmer originally considered mostly $C^{\infty}$ actions.   Most evidence for the  conjecture including existing results requires the action to be at least  $C^1$ but the conjecture might be true for actions by homeomorphisms, see particularly \cite{MR2807834, MR3150210} for discussion and evidence in this regularity.
\end{HistoricalRemarks}

The group $\Sl(n,\R)$ is the standard split simple Lie group with restricted root system of type $A_n$. We denote by  $\Sp(2n,\mathbb R)$ the  group of real symplectic $2n\times 2n$ matrices, the standard split simple Lie group of rank $n$ with restricted root system of type $C_n$.

\begin{theorem}
\label{thm:sp(n,R)}
Conjecture \ref{conjecture:zimmer} holds for cocompact lattices in $\Sp(2n,\mathbb R)$ for $n\ge  2$.  In particular if $M$ is a compact manifold with  $\dim(M)< 2n-1$ and $\Gamma < \Sp(2n,\R)$ is a cocompact lattice then any  homomorphism $\alpha\colon \Gamma \rightarrow \Diff^2(M)$ has finite image. In addition, if $\dim(M)=2n-1$ and $\omega$ is a volume form on $M$  then any  homomorphism $\alpha\colon \Gamma \rightarrow \Diff^2(M, \omega)$ has finite image.
\end{theorem}

The fact that all actions in Theorem \ref{slnr} and \ref{thm:sp(n,R)} factor through finite quotients
follows from the existence of an invariant Riemannian metric and the fact that, for these cases, $v(G)+1=n(G) \le d(G)$
where $v(\Sl(n,\R)) = n-1$ and $v(\Sp(n,\R)) = 2n-1$.  See  Section \ref{section:finalarguments} for full discussion.

The remaining  split simple classical Lie groups are $\So(n,n)$ and $\So(n,n+1)$.  Note that $\So(2,2)$  is not simple and we omit below the higher-rank simple groups  $\So(2,3)$ and $ \So(3,3)$ as their identity components are   double covered by $\Sp(4,\R)$  and $\Sl(4,\R)$, respectively.
For $G= \SO(n,n)$ with $n\ge 4$, we have
	$$n(G) = 2n ,\quad  d(G) = 2n-1, \text{ and } \quad  v(G) = 2n-2\phantom{.}$$
	and similarly for $G= \SO(n,n+1)$ with $n\ge 3$, we have
	$$n(G) = 2n+1, \quad d(G = 2n,  \text{ and } \quad v(G) = 2n-1 .$$

\begin{theorem}
\label{thm:orthogonalgroups}
The non-volume-preserving case of Conjecture \ref{conjecture:zimmer} holds for cocompact lattices $\Gamma$ in $\So(n,n)$ with $n\ge 4$ and for $\So(n,n+1)$ with $n\ge 3$; the volume-preserving case holds up to dimension 1 less than conjectured.  

More precisely, let $M$ be a compact connected manifold and $\omega$ a volume form on $M$.  
\begin{enumerate}
\item If \, $\Gamma<\So(n,n)$ is a cocompact lattice and  $\dim(M)< 2n-2$ then any  homomorphism $\alpha\colon \Gamma \rightarrow \Diff^2(M)$ has finite image.  If $\dim(M) = 2n-2$ 
then any  homomorphism $\alpha\colon \Gamma \rightarrow \Diff^2(M,\omega)$ has finite image. 
\item If \, $\Gamma<\So(n,n+1)$ is a cocompact lattice and  $\dim(M)< 2n-1$ then any  homomorphism $\alpha\colon \Gamma \rightarrow \Diff^2(M)$ has finite image.  If $\dim(M) = 2n -1$  
then any  homomorphism $\alpha\colon \Gamma \rightarrow \Diff^2(M,\omega)$ has finite image.  
\end{enumerate}
\end{theorem}

Again, the finiteness of the action  follows from  Theorem \ref{theorem:main} below and a computation of the value of  $d(G)$.

From  Conjecture \ref{conjecture:zimmer} for split orthogonal groups, one expects that in dimension $n(G)-1= d(g)=v(g) +1$ all volume-preserving actions necessarily preserve  a Riemannian metric.   In this case,  Margulis' superrigidity theorem would imply the action is finite unless the manifold is the $(n(G)-1)$-dimensional sphere or projective space.  While the techniques of this paper impose certain restrictions on   volume-preserving  actions in dimension $n(G)-1$,   it seems additional ideas are needed to obtain the conjectured result in dimension $n(G)-1$.

We remark that the conclusions of Theorems \ref{slnr}, \ref{thm:sp(n,R)}, and \ref{thm:orthogonalgroups} continue to hold
for actions of cocompact lattices in connected Lie groups   isogenous to the groups in the theorems.  That is, if $G$ is a connected Lie group with finite center whose Lie algebra is isomorphic to the Lie algebra of a group in Theorems \ref{slnr}, \ref{thm:sp(n,R)}, or \ref{thm:orthogonalgroups}, then the conclusion  of the corresponding theorem continues to hold for cocompact lattices in $G$.

Combined with the main results of \cite{MR2219247} and \cite{AWBFRHZW-latticemeasure} we obtain the following theorem for actions of  lattices on surfaces.

 \begin{theorem}[{\cite[Corollary 1.7]{MR2219247} + \cite[Theorem 1.6]{AWBFRHZW-latticemeasure} + Theorem \ref{theorem:main}}]\label{thm:easycor2}
Let $S$ be a closed, oriented surface of genus at least 1.  Let $G$ be a connected simple Lie group with finite center and real-rank at least $2$ and assume the  restricted root system of the Lie algebra of $G$ is not of type $A_2$.    Let $\Gamma\subset G$ be a lattice.  Then any  homomorphism $\alpha\colon \Gamma \rightarrow \Diff^2(S)$ has finite image.
\end{theorem}

  Note that the hypothesis that the restricted root system of $G$ is not of type $A_2$ ensures the number $r(G)$ defined in Section \ref{subsection:algebra} below is at least $3$.  Up to isogeny, the  three simple Lie groups of type $A_2$ are $\Sl(3,k)$ where $k= \R, \C$, or $\mathbb H$.  We remark that the conclusion of Theorem \ref{thm:easycor2} is expected to hold for lattices in $\Sl(3, \C)$ and $\Sl(3, \mathbb H)$, and  for lattices in $\Sl(3,\R)$ assuming that  $S$ is not the 2-sphere.

We defer the statement of our main theorem, Theorem \ref{theorem:main}, which includes partial results for non-split  and exceptional Lie groups, until we have made some requisite definitions.
For non-split groups, our main theorem does not   recover the full conjecture but does imply finiteness of actions in a dimension that grows linearly with the rank.

\subsection{Outline of the proof}
We will illustrate the main ideas of the proof of Theorem \ref{slnr} by considering the case where $\Gamma \subset G = \Sl(n, \R)$ is a cocompact lattice acting on a closed manifold $M$ and $\text{dim}(M) < n-1$. In this case, if the action preserves a measure $\mu$, Zimmer's cocycle superrigidity theorem implies that the derivative cocycle is measurably cohomologous to a cocycle taking values in a compact subgroup or, equivalently, that the action preserves a measurable Riemannian metric \cite{MR595205}. This implies, in particular, that all Lyapunov exponents for all elements of $\Gamma$ are zero.   As remarked above, the conjecture would follow from promoting the invariant measurable metric  to a smooth invariant metric.



It was observed by Zimmer that  conjecture would follow from the existence of an invariant Riemannian metric of quite low regularity.  Indeed, in the case of volume preserving actions, Zimmer observed  it sufficed to find a metric that was bounded above and below in comparison to a background smooth metric; that is, it suffices to find an invariant  $L^{\infty}$ metric.
Very early on, Zimmer also observed that one might get better regularity by noting that the metric was invariant, so its growth along orbits was controlled by the derivative cocycle. Using this he could show that the metric was, in a sense, in $L^{\epsilon}$, for  very small values of $\epsilon >0$ \cite{MR743815}.  A more sophisticated, non-linear, attempt to average metrics in order to produce invariant smooth metrics was proposed by the second author in \cite[Section 4.6.2]{F11}.  Both of these attempts fail to produce good results because even with a measurable (or even slightly more regular) invariant metric, the only a priori bound on growth of derivatives along orbits is exponential.

The first   step in the proof of Theorem \ref{slnr} is to show  that any action $\alpha\colon \Gamma\to \diff^2(M)$ for $\Gamma$ and $M$ as in Theorem \ref{slnr} has \emph{uniform subexponential growth of derivatives}:  for every $\e>0$, there is $C_{\e}$ such  that  $$\|D\alpha(\gamma)\| \leq C_{\e}e^{\e \text{l}(\gamma)}$$
 where $\|D\alpha(\gamma)\| = \text{max}_{x \in M} \|D_x\alpha(\gamma)\|$ denotes the norm of the derivative and $\text{l}(\cdot)$ denotes the word-length with respect to some choice of finite generating set for $\Gamma$.

To illustrate how we establish uniform  subexponential growth of derivatives, consider a  more elementary fact from classical smooth dynamics:  a  diffeomorphism $f\colon M\to M$ of a compact manifold $M$  has uniform subexponential growth of derivatives if and only if all   Lyapunov exponents of $f$ are zero with respect to any $f$-invariant probability measure.   Clearly, uniform  subexponential growth of derivatives implies  that all Lyapunov exponents vanish for any measure.   To prove the converse, assume that for some fixed $\epsilon >0$ there are $x_n$ and $N_n \rightarrow \infty$ so that  $\|D_{x_n}f^{N_n}\| \geq e^{\e {N_n}}$; then any accumulation point $\mu$ of the  sequence of measures $\mu_n := \frac{1}{N_n} \sum_{i=1}^{N_n} f^i *\delta_{x_n}$  will be a measure $\mu$ whose average top Lyapunov exponent  (see discussion in Section \ref{subsection:lyapunov} and  \eqref{eq:toplyap}  below) is positive.

To implement the above idea  in the context of $\Gamma$-actions rather than $\mathbb{Z}$-actions, in Section \ref{subsection:suspension} we   induced from the $\Gamma$-action on $M$ to a $G$-action on an auxiliary manifold 
  $M^\alpha$.  
  This space has the structure of an $M$-bundle over $G/\Gamma$.  For $A\subset \Sl(n,\R)$ the subgroup of positive diagonal matrices (that is, a maximal split Cartan subgroup), the failure of the action $\alpha$ to have uniform subexponential growth of derivatives implies the existence of an element $s\in A$ and an $s$-invariant probability measure  $\mu$ on $M^\alpha$ with a positive Lyapunov exponent for the fiberwise derivative cocycle.    The key new idea  is to construct from $\mu$ a $G$-invariant measure $\mu'$ on $M^\alpha$ such that the fiberwise derivative cocycle continues to have a positive Lyapunov exponent for some $s'\in A$.  This yields a contradiction with Zimmer's cocycle superrigidity theorem as there are no non-trivial linear representations in dimension less then $n$.  We thus obtain the uniform subexponential growth of derivatives for the action $\alpha$.

To construct a $G$-invariant measure $\mu'$, starting with our $s$-invariant measure $\mu$ we build a sequence of measures by averaging: given a measure $\mu$ that has a positive fiberwise Lyapunov exponent for some $s\in A$, by averaging $\mu$ along $A$ or a unipotent subgroup  commuting with $s$ we obtain a new measure $\mu'$ with better invariance properties   and has positive fiberwise exponent for some $s'\in A$.  There is some similarity here to Margulis' original proof of the superrigidity theorem using Oseledec's theorem where it is used  (see  \cite{MR1090825}) that higher-rank semisimple Lie groups can be generated by centralizers of certain  elements of the diagonal subgroup.

While we cannot average directly to obtain  a $G$-invariant measure on $M^\alpha$, we may average so as to obtain an $A$-invariant measure on $M^\alpha$ whose projection to $G/\Gamma$ is the Haar measure and has positive fiberwise exponent for some $s'\in A$ .  This step  requires a careful choice of subgroups  over which to average and employs Ratner's theorem on measures invariant under unipotent subgroups and an improvement due to Shah concerning averages of measures along unipotent subgroups.  As the general averaging argument requires  understanding the combinatorics of   root systems,   we explain this step for the special case of $\Sl(n, \R)$ in Section \ref{sec:SL}.

To show such a  measure is, in fact, $G$-invariant,  we use a result  (Proposition \ref{thm:nonresonantimpliesinvariant} below) from the work of the first author with Rodriguez Hertz and Wang where its shown that---under the same dimension bounds as in Theorem \ref{slnr}---any $P$-invariant measure on $M^\alpha $ is, in fact, $G$-invariant \cite{AWBFRHZW-latticemeasure}. Here $P$ denotes the group of upper triangular matrices.
As $P$ contains $A$ and as any $P$-invariant measure on $G/\Gamma$ is necessarily Haar, we are in a slightly more general setting than considered in \cite{AWBFRHZW-latticemeasure}.  The key idea in the proof  in \cite{AWBFRHZW-latticemeasure} of Proposition \ref{thm:nonresonantimpliesinvariant}   is to relate the Haar-entropy  of elements of the $A$-action on $G/\Gamma$ with    the $\mu$-entropy of elements of the $A$-action on $M^\alpha$.  For the Haar measure on $G/\Gamma$, the   entropy of elements of $A$ is   computed in terms of the roots of $G$.   Moreover,  the contribution from the fiber to the $\mu$-entropy of elements of  the $A$-action is constrained by the dimension assumption.  Many key ergodic theoretic notions for these argument are developed in \cite{AWB-GLY-P1, AWB-GLY-P2,AWB-GLY-P3}.  

Both the main result in \cite{AWBFRHZW-latticemeasure} and our use of their techniques here employ the philosophy that ``non-resonance implies invariance." This philosophy was introduced by the same authors in their study of global rigidity of Anosov actions of higher-rank lattices in  \cite{BRHW-Anosov}.  Given a $G$-action and an $A$-invariant (or equivariant) object $O$, such as a measure or a semiconjugacy to a linear action, one may try to associate to $O$ a class of linear functionals $\mathcal O$. In the case of an $A$-invariant measure, the functionals are the Lyapunov exponents; in the case of a conjugacy to a linear action, the functionals are the weights of the representation corresponding to the linear action.
The philosophy,  implemented in both  \cite{BRHW-Anosov} and  \cite{AWBFRHZW-latticemeasure}, is that, given any root $\beta$ of $G$ that is not positively proportional to an element of $\mathcal O$, the object $O$ will automatically be invariant (or equivariant) under the unipotent subgroup associated to $\beta$ (or to $\beta$ and $2\beta$).   If one can find enough such non-resonant roots, the object $O$ is automatically $G$-invariant (or $G$-equivariant).

The second   step in the proof of Theorem \ref{slnr} is to use  strong property (T) introduced by V.~Lafforgue and uniform subexponential growth of derivatives to produce an invariant metric for the action.   Strong property (T) was introduced by Lafforgue who proved that all simple Lie groups containing $\Sl(3,\mathbb R)$ and their cocompact lattices have strong property $(T)$ with respect to Hilbert spaces. The precise results we use here are an extension of Lafforgue's due to de Laat and de la Salle \cite{MR2423763,MR3407190}.

We formulate a special case of the results of \cite{MR2423763,MR3407190} below.  Given a Hilbert or Banach space $\mathcal{H}$, let $B(\mathcal{H})$ denote the bounded operators on $\mathcal{H}$.
\begin{theorem}[\cite{MR3407190}]\label{lafforgue1} Let $\mathcal{H}$ be a Hilbert space and let $\Gamma$ be as in Theorem \ref{slnr}.  There exists $\e>0$, such that for any representation $\pi\colon \Gamma \to B(\mathcal{H})$, if there exists $C_\e>0$ such that
\begin{align*}
\|\pi(g)\| \leq C_{\e}e^{\e \text{l}(\gamma)},
\end{align*}
then there exists a sequence of averaging operators $p_n =  p(\mu_n)$ in $B(\mathcal{H})$, defined by probability measures $\mu_n$ on $\Gamma$ supported in the ball of radius $n$, such that for any vector $v\in \mathcal{H}$, the sequence $v_n = p_n(v) \in \mathcal{H}$ converges to a $\Gamma$-invariant vector $v^{*}$. Moreover the convergence is exponentially fast: there exists  $0< \lambda <1$  (independent of $\pi$) and a $C$ so that  $\|v_n-v_*\| \leq C \lambda^n \|v\|$.  
\end{theorem}

In the case of $C^\infty$ actions, we may apply this theorem to the Sobolev space of sections of the bundle of symmetric $2$-tensors on $M$ (which contains the space of Riemannian metrics as a subset).  As the  uniform subexponential growth of derivatives implies subexponential growth of derivatives of higher order (see  Lemma \ref{lemma:estimateoncompostion} below), we  verify the slow norm growth required in Theorem \ref{lafforgue1}. Starting from an initial symmetric 2-tensor field $g$ which is a Riemannian metric, we  obtain from Theorem \ref{lafforgue1} a non-negative,  $\Gamma$-invariant,  symmetric 2-tensor field on $M$.  To verify that the tensor is in fact a metric (that is, to verify that the 2-tensor is  non-degenerate) we use that the norms decay at a subexponential rate under the averaging operator while the convergence to the limit is exponentially fast.

We remark that a somewhat similar use of subexponential growth of derivatives along a central foliations also occurs in the work of the second author with Kalinin and Spatzier on rigidity for Anosov actions of abelian groups \cite{MR2983009}. In that work, subexponential growth is verified from the existence of  a H\"{o}lder conjugacy and is used in conjunction with exponential decay of matrix coefficients for abelian groups.  These ideas are also applied in the work of Rodriguez-Hertz and Wang \cite{MR3260859}.

\subsection*{Acknowledgements}
We thank Yves Benoist, Alex Eskin, Federico Rodriguez Hertz, Mikael de la Salle, Amie
Wilkinson, Dave Witte Morris, and Bob Zimmer for useful conversations. We thank the referees, Homin Lee and Federico Vigolo for many useful comments on earlier versions of this paper.  Particular thanks are due to an anonymous referee for suggesting the proof of Proposition \ref{prop:LieFacts2}, to Dave Witte Morris for suggesting the proof of Lemma \ref{lem:dave} and to Bob Zimmer for his enthusiasm and encouragement.

\section{Main  theorem and proof of results from  introduction}
Our main theorem, Theorem \ref{theorem:main} below, gives a partial solution to Zimmer's conjecture for actions of cocompact lattices in any semisimple Lie group all of whose non-compact, almost-simple factors are of  higher rank.  
Results  stated in the introduction  follow from Theorem \ref{theorem:main} 
 and Margulis' superrigidity theorem  as explained below in Section \ref{section:finalarguments}.

\subsection{Main theorem}
To state our main theorem, given a semisimple Lie group we associated an integer  $r(G)$ similar to $v(G)$ in Conjecture \ref{conjecture:zimmer}.
For $\R$-split Lie groups $G$, we always have $r(G) = v(G)$.  More generally, we have $r(G) = v(G')$, where $G'$ is a maximal $\R$-split simple subgroup of~$G$.  
An alternative definition of $r(G)$ in terms of root data  is given below in Definition \ref{def:rescodim}.
\label{subsection:genresults}
\begin{theorem}\label{theorem:main}
Let $G$  be a connected, semisimple real Lie group with finite center, all of whose non-compact,  almost-simple factors have real-rank at least $2$.  Let $\Gamma\subset G$ be a cocompact lattice and for $k \ge 2$,  let $\alpha \colon \Gamma\to \diff^k(M)$ be an action.  Suppose that either
\begin{enumerate}
\item $\dim(M)<r(G)$, or
\item $\dim(M)= r(G)$ and $\alpha$ preserves a smooth volume.
\end{enumerate}
Then $\alpha(\Gamma)$ preserves a Riemannian metric which is $C^{k-1-\delta}$ for all $\delta>0$.
\end{theorem}

  Theorem \ref{theorem:main} gives a partial solution to Zimmer's conjecture  for cocompact lattices in any higher-rank simple Lie group $G$.  In particular, the number $r(G)$ provides a    critical dimension---which grows linearly in the rank of $G$---for which the conclusion of Zimmer's  conjecture holds.  Moreover, the number $r(G)$ gives the optimal result  for non-volume-preserving actions when  $G$ is a split real form.

For non-split simple Lie groups, our critical dimension falls below the conjectured result.
In particular, while we   recover the complete conjecture as stated in Conjecture \ref{conjecture:zimmer}  for cocompact lattices in $\Sl(n,\R)$ with $n>2$, for lattices in $\Sl(n,\C)$ and $\Sl(n,\mathbb{H})$ our critical dimension $r(G)$ is, respectively,  one half and one quarter of the conjectured critical value.
For  lattices in $\So(n,m)$ we obtain the conjectured result in the split case where $m=n$ or $m=n+1$.
However, for fixed $n$ our critical dimension  $r(G)$ for $G= \So(n,m)$, $m> n$, is constant in $m$ and thus the defect between the critical dimension in   Theorem \ref{theorem:main} and the  conjectured critical dimension  becomes arbitrarily large   as $m\to \infty$.

The  obstruction to improving our results for non-split simple Lie groups is to improve the results of \cite{AWBFRHZW-latticemeasure},
particularly the result quoted below in Proposition \ref{thm:nonresonantimpliesinvariant}.  In particular, the method of proof of  Proposition \ref{prop:Ginvariantmeasure} below can not distinguish between actions of lattices in two groups with the same restricted root system.

\begin{remark}
\label{rmk:nocompacts}
In   Theorem \ref{theorem:main}  above, by restricting to a finite-index subgroup of  $\Gamma$  it is with no loss of generality to assume the group $G$ has no compact factors and is center free.  Indeed, $G$ is an almost direct product $G= KL$ where $K$ is the largest compact normal subgroup of $G$ and $L$ has no compact normal subgroups of positive dimension.   Since compact groups are linear, the image of $\Gamma$ in $G/L$ has a torsion-free subgroup of finite index.  Then there is a finite-index subgroup $\Gamma'$ of $\Gamma$ such that $\Gamma'\cap K$ is the identity.  Then, the map $G\to G/K$ restricts to an injection on $\Gamma'$; thus an action of the subgroup $\Gamma$  of $G$ induces and action of the subgroup $\Gamma'$ of $G/K$.
\end{remark}

In the remainder of the paper, we will assume $G$ has no compact factors to simplify some algebraic arguments.

\subsection{Proof of Theorem \ref{theorem:main}}
We prove  Theorem \ref{theorem:main}  in  two steps.

Let $\Gamma$ be a finitely generated group.  Let $\wl\colon \Gamma\to \N$ denote the word-length function relative to some fixed  finite symmetric set of generators.
Let $\alpha \colon \Gamma \to \diff^1(M)$ be an  action of $\Gamma$  on a compact manifold $M$ by $C^1$ diffeomorphisms.  
We say the action $\alpha$ has \emph{uniform subexponential growth of derivatives} if for all $\epsilon>0$ there is a $C_{\epsilon}$ such that for all $\gamma \in \Gamma$ we have
$$\|D\alpha(\gamma)\|\le C_{\epsilon}e^{\epsilon \wl(\gamma)}$$
where $\|D\alpha(\gamma)\|= \sup_{x \in M} \|D_x\alpha(\gamma)\|$.

To prove Theorem \ref{theorem:main} we first establish uniform subexponential growth of derivatives  for actions of cocompact lattices in  the low-dimensional settings consider above.
\begin{theorem}\label{prop:USEGOD}
Let $G$  be a connected, semisimple Lie group with finite center.    
Let $\Gamma\subset G$ be a cocompact lattice and let $\alpha \colon \Gamma\to \diff^{1+\beta}(M)$ be an action for $\beta>0$.  Suppose that either
\begin{enumerate}
\item $\dim(M)<r(G)$, or
\item $\dim(M)= r(G)$ and $\alpha$ preserves a smooth volume.
\end{enumerate}
Then $\alpha$ has uniform subexponential growth of derivatives.
\end{theorem}
When $G$ is rank-1 or has rank-1 factors we have  $r(G)=1$.  In this case,  Theorem \ref{prop:USEGOD} is  trivial if $\dim (M) <r(G)$ and  is nearly as  trivial if $\dim(M)= r(G)$ and $\alpha$ preserves a smooth volume since any group of diffeomorphisms preserving a smooth volume form on the circle is smoothly conjugate to a group of isometries.

Having established Theorem \ref{prop:USEGOD}, the second step in the proof of Theorem \ref{theorem:main} is to show that for a group with strong property $(T)$,  any action with subexponential growth of derivatives preserves a smooth Riemannian metric.

\begin{theorem}
\label{thm:strongTplussubexponential}
Let $\Gamma$ be a finitely generated group, $M$ a compact manifold,  and $\alpha \colon \Gamma \rightarrow \Diff^k(M)$ an action on $M$ by $C^k$ diffeomorphisms for $k \geq 2.$   If $\Gamma$ has strong property (T) and if  $\alpha$ has  uniform subexponential growth of derivatives then $\alpha$ preserves a Riemannian metric which is $C^{k-1-\delta}$ for all $\delta >0$.
\end{theorem}

Theorem \ref{theorem:main} is an immediate consequence of Theorems \ref{thm:strongTplussubexponential} and \ref{prop:USEGOD}.


Note that Theorem \ref{theorem:main}  implies Conjecture \ref{conjecture:zimmer} for  non-volume-preserving actions of cocompact lattices in all split simple Lie groups.   Moreover, as the minimal  non-trivial linear representations of $\mathfrak{sl}(n,\R)$ and $\mathfrak{sp}(2n,\R)$ occur in dimensions $n$ and $2n$, respectively, Theorem \ref{theorem:main} implies the volume-preserving case of  Conjecture \ref{conjecture:zimmer} for lattices in (groups isogenous to) $\Sl(n,\R)$ and $\Sp(2n, \R)$.  For the split orthogonal groups, the minimal linear representations occur in dimensions $2n = r(\lieg) +2$ for $\lieg = \mathfrak{so}(n,n)$ and $2n+1 = r(\lieg) +2$ for $\lieg = \mathfrak{so}(n,n+1)$ and thus
we are unable to  recover the full  conjecture   for volume-preserving actions from Theorem \ref{theorem:main}.

\subsection{From metrics to compact Lie groups and finite actions}
\label{section:finalarguments}
To complete the proofs of the results from the introduction we recall that the isometry group of the metric guaranteed by Theorem \ref{theorem:main} is a compact Lie group whose dimension is bounded from above; we then apply Margulis's superrigidity theorem with compact codomains to conclude the image $\alpha(\Gamma)$ is finite.  All arguments in this subsection are well known to experts and we include them for completeness.

Let $M$ be equipped with the metric guaranteed by Theorem \ref{theorem:main}.  We claim the isometry group of $M$
is a compact Lie group.  When the metric
is at least $C^2$ this is an immediate consequence of the Myers-Steenrod Theorem and the fact that $\Isom(M)$ embeds in the bundle
of orthogonal frames over $M$ which is an $O(\dim(M))$ bundle \cite{MR1503467, MR1336823}.  When the metric is not $C^2$ an additional argument is needed to show
that isometries are at least $C^1$.  Recently Taylor proved that isometries of an $\alpha$-H\"{o}lder Riemannian metric are $C^{1+\alpha}$ \cite{MR2204038}. See also related work in \cite{MR283727,MR2216252}.  Given Taylor's result, we again have an embedding of $\Isom(M)$ into the bundle of orthogonal frames and so $\Isom(M)$ is a compact Lie group.  One can also argue instead by viewing $M$ as a compact metric space whose isometry group is compact and use the resolution of Lipschitz case of  the Hilbert-Smith conjecture by Repov{\u{s}} and \u{S}{\u{c}}epin to see $\Isom(M)$ has no small subgroups and is therefore   a Lie group \cite{MR1464908}.  Isometries of a H\"{o}lder Riemannian metric, or even an $L^{\infty}$ Riemannian metric, are easily seen to be Lipschitz maps.



We now prove finiteness of the image $\alpha(\Gamma)$ in any theorem from the introduction.
We assume that $\alpha\colon  \Gamma \rightarrow \Isom(M)$ and show that if $\alpha(\Gamma)$ is infinite then $\dim(M) \geq d(G)$.  Let $L= \overline {\alpha (\Gamma)}$ be the closure of $\alpha(\Gamma)$ in $\Isom(M)$. Passing to a finite index subgroup of $\Gamma$ one can assume $L$ is connected. By the structure theory of compact Lie groups $L$ is an almost direct product $L=K_1 \times  \cdots \times K_r$.  Using that compact groups are real algebraic and applying Margulis's superrigidity and arithmeticity theorems we will see that each $K_i$ is a compact real form of a simple factor of the complexification of $G$.   First, since the abelianization of $\Gamma$ is trivial, all factors of $L$ are simple.  To prove all remaining assertions, we need only consider a single factor $K=K_i$.  Let $H$ be the complexifation of $K$.  Arguing as in \cite[Lemma 6.1.6]{MR0776417}, we can see that the trace of $\Ad_{\mathfrak{h}}(\alpha(\gamma))$ is real algebraic for every $\gamma \in \Gamma$.   We can then apply \cite[Lemma 6.1.7]{MR0776417} to find an embedding of $K$ into $GL(n,\mathbb{C})$ such that each $\alpha(\gamma)$ has algebraic entries.  Since $\Gamma$ is finitely generated, it follows that there is a number field $k$ such that each $\alpha(\gamma)$ has entries in $k$ and $K$ is defined over $k$. Applying superrigidity with $p$-adic targets, we see that $\Gamma$ has a finite-index subgroup for which every $\alpha(\gamma)$ has entries lying in the integer points $\mathcal{O}_k$ of $k$.  Applying restriction of scalars, \cite[Proposition 6.1.3]{MR0776417}, we see that either $\alpha(\Gamma)$ is contained in the integer points of a compact group and  is thus finite, or there is a field automorphism $\sigma$ of $k$ over $\mathbb{Q}$ such that $\sigma(\alpha(\Gamma))$ is Zariski dense and unbounded in a non-compact simple group $G'$. Applying Margulis' superrigidity theorem again,  $G'$ is locally isomorphic to a factor of $G$.   Since $G'$  is locally isomorphic to a factor of $G$, restriction of scalars implies that $K$ is a compact real form of a simple factor of the complexification of $G$.  Furthermore since $K< \Isom(M)$ is non-trivial, there is a closed $K$ orbit in $M$ of the form $K \cdot x =K/C$ for some proper subgroup $C$.  This then forces $\dim(M) \geq \dim(K/C) \geq d(G)$.

To complete the proofs of the results in the introduction, one computes the number  $d(G)$ appearing in Conjecture \ref{conjecture:zimmer}, the minimal dimension of $K/C$ where $K$ is a compact real form $K$ of the classical Lie  group $G$ and $C$ is a proper closed subgroup.  In all cases considered in the introduction,  $d(G) > \dim M$ and finiteness of the action follows.   


\section{Background and facts from Lie theory}
We recall some facts and definitions from the structure theory of real Lie groups as well as some notation that will be used in the sequel.  A standard reference is \cite{MR1920389}.  For the reader interested only in actions of cocompact lattices in $\Sl(n,\R)$, we recommend skipping this section on the first read.
\subsection{Structure theory   of Lie groups} 
Let $G$ be a connected,  semisimple Lie group with finite center.
As usual, write $\lieg$ for the Lie algebra of $G$.  Fix once and for all a Cartan involution  $\theta$  of $\lieg$ and write $\liek$ and $\liep$, respectively,  for the $+1$ and $-1$ eigenspaces of $\theta$.  Denote by $\liea$  a maximal abelian subalgebra of $\liep$ and by $\liec$ the centralizer of $\liea$ in $\liek$.  We let $\Sigma$ denote the set of
restricted roots of $\lieg$ with respect to $\liea$. Note that the elements of $\Sigma$ are real linear functionals on $\liea$.
Recall that $\dim_\R(\liea)$ is the real-rank of $G$.  We fix $\liea$ for the remainder.  

Recall that a \emph{base} (or a collection of \emph{simple roots}) for $\Sigma$ is a subset $\Pi\subset \Sigma$ that is a basis for the vector space $\liea^*$ and such that every non-zero root $\beta\in \Sigma$ is either a positive or a negative integer combination of elements of $\Pi$.  For a choice of $\Pi$, elements $\beta\in \Pi$ are called \emph{simple} (positive) roots.  Relative to a choice of base $\Pi$, let $\Sigma_+\subset \Sigma$ be the collection of positive roots
and let $\Sigma_-$ be the corresponding set of negative roots.  For $\beta\in \Sigma$ write $\lieg^\beta$ for the associated root space.  Then $\lien =\bigoplus_{\beta\in \Sigma_+} \lieg^\beta$ is a nilpotent subalgebra.
A  subalgebra  $\lieq$ of $\lieg$ is said to be a \emph{standard parabolic subalgebra} or simply \emph{parabolic} (relative to  the choice of $\theta$ and $\Pi$) if 
 $ \liec\oplus\liea\oplus \lien\subset \lieq $ where $\lien$ is defined relative to $\Pi$.
We have that the standard parabolic subalgebras of $\lieg$ are parametrized by exclusion of simple (negative) roots: for any sub-collection $\Pi'\subset\Pi$ let
\begin{equation}\label{eq:parabolic}\lieq_{\Pi'}= \liec\oplus \liea \oplus \bigoplus_{\beta\in \Sigma_+ \cup \mathrm{Span}(-\Pi')} \lieg^\beta.\end{equation}
Then $\lieq_{\Pi'}$ is a Lie subalgebra  of $\lieg$ and all standard parabolic subalgebras of $\lieg$ are of the form $\lieq_{\Pi'}$ for some $\Pi'\subset \Pi$  \cite[Proposition 7.76]{MR1920389}.

Let $A,N,$ and $ K$ be the {analytic subgroups} of $G$ corresponding to $\liea, \lien$ and $\liek$.
Then $G = KAN$ is the corresponding Iwasawa decomposition of $G$.  As $G$ has finite center, $K$ is compact.      Note that the Lie   exponential $\exp\colon \lieg\to G$ restricts to  diffeomorphisms between $\liea$ and $A$ and $\lien$ and $N$.  Fixing a basis for $\liea$, we often identify $A= \exp (\liea) = \R^d$.  Via this identification we    extend linear functionals on $\liea$ (in particular, the restricted roots of $\lieg$) to functionals on $A$.  
Write $C$ for the centralizer of $\liea$ in $K$ and recall that $\mathfrak{c}$ is the Lie algebra of $C$.
Then $P = CAN$ is the \emph{standard minimal parabolic subgroup}.
 Since $C$ is compact, it follows that $P$ is amenable.
 A \emph{standard parabolic subgroup} (relative to the choice of $\theta$ and $\Pi$ above) is any closed subgroup  $Q\subset G$ 
 containing $P$. 
 The Lie algebra of any standard parabolic subgroup $Q$ is a standard parabolic subalgebra and the correspondence between standard parabolic subgroups and subalgebras is 1-1.

We say two restricted roots $\beta,\hat \beta\in \Sigma$ are \emph{positively proportional} if there is some $c>0$ with
$$\hat\beta = c\beta.$$
Note that $c$ takes values only  in $\{\frac 1 2, 1, 2\}$ and this occurs only if the root system $\Sigma$ has a factor of type $BC_\ell$.
Let $\hat\Sigma$ denote the set of \emph{coarse restricted roots}; that is, $\hat\Sigma$ denotes the collection of positively proportional equivalence classes $[\beta]$ in $\Sigma$.
Note that for $[\beta] \in \hat \Sigma$, $\lieg^{[\beta]} := \bigoplus _{\beta' \in [\beta]} \lieg ^{\beta'}$ is a nilpotent subalgebra      
 and the Lie exponential  map  restricts to a  diffeomorphism between $\lieg^{[\beta]}$ and  the corresponding analytic subgroup which we denote by $G^{[\beta]}$.

Let $\lieq$ denote a standard parabolic subalgebra of $\lieg$.
Observe that if $\lieg^\beta\cap \lieq\neq 0$ for some $\beta\in \Sigma$ then, from the structure of parabolic subalgebras, $\lieg^{[\beta]}\subset \lieq$ where $[\beta]\in \hat \Sigma$ is the coarse restricted root containing $\beta$.
A     proper subalgebra $\lieh$ of $\lieg$ is \emph{maximal} if there is no subalgebra $\lieh'$ with $\lieh\subsetneq\lieh'\subsetneq\lieg$.  Note that maximal standard parabolic subalgebras are of the form $\lieq_{\Pi\sm \{\beta\}}$ for some $ \beta\in \Pi$.

\subsection{Resonant codimension and related lemmas}
\label{subsection:algebra}

We say a Lie algebra is {\em saturated by coarse roots spaces} if its intersection with a coarse root space
is either trivial or the entire coarse root space.
Consider a Lie subalgebra $\lieh\subset \lieg$ that is saturated by coarse root spaces.  For such a subalgebra define the \emph{resonant codimension}, $\bar r (\lieh)$, of $\lieh$ to be the cardinality of the set $$\{[\beta] \in \hat \Sigma \mid \lieg^{[\beta]}\not\subset \lieh\}.$$    For a subgroup $H\subset G$ whose Lie algebra is saturated by coarse root spaces, we will also refer to the resonant codimension of the group $H$.

Note that standard parabolic subalgebras $\lieq$ are automatically saturated by coarse root spaces whence the {resonant codimension}  is defined for all standard parabolic subalgebras.
 As in \cite{AWBFRHZW-latticemeasure}), given a (semi)simple Lie algebra $\lieg$ as above we define a combinatorial number $r(\lieg)$.  As the number depends only on the Lie algebra $\lieg$, we use both the notation  $r(G)$ and $r(\lieg)$ interchangeably.

\begin{definition}\label{def:rescodim}The \emph{minimal resonant codimension} of $\lieg$ or $G$, denoted by $r(\lieg)$ or $r(G)$, is defined to be the minimal value of the resonant codimension $\bar r (\lieq)$  of $\lieq$ as $\lieq$ varies over all (maximal)   proper  parabolic subalgebras of $\lieg$.
\end{definition}

\begin{remark}\label{rem:ooppllkk} In the case that the Lie algebra $\lieg$ of $G$ is a split real form, the minimal resonant codimension $r(\lieg)$ coincides with minimal codimension of all maximal parabolic subalgebras.  In general, we have $r(\lieg) \leq v(G)$.  That this definition of $r(G)$ agrees with the one given immediately before Theorem \ref{theorem:main} follows from \cite[Theorem 7.2]{MR0207712}.

In the case that $\lieg$ is semisimple then $r(\lieg)$ is the minimal value of $r(\lieg')$ as $\lieg'$ varies over all non-compact simple ideals of $\lieg$.  In particular, if $\lieg$ has rank-1 factors then $r(\lieg)=1$.
\end{remark}

\begin{example}We compute $r(\lieg)$ for a number of classical real simple Lie algebras.
Note that it follows  from  definition that the minimal resonant codimension  depends only on the restricted root system of $\lieg$ and  not on the Lie algebra $\lieg$.
\begin{description}
\item [Type $A_n$] $r(\lieg)= n$.  This includes $\mathfrak{sl}(n+1,\R),$ $ \mathfrak{sl}(n+1,\mathbb{C}), $ $\mathfrak{sl}(n+1,\mathbb H)$.
\item [Type $B_n$, $C_n$, and $(BC)_n$] $r(\lieg)= 2n-1$. This includes  $\mathfrak{sp}(2n,\R)$, $\mathfrak{so}(n,m)$ for $n<m$, and
$\mathfrak{su}(n,m)$ and
$\mathfrak{sp}(n,m)$ for $n\le m$.
\item [Type $D_n$, $n\ge 4$] $r(\lieg)= 2n-2$.  This includes $\mathfrak{so}(n,n)$ for $n\ge 4$
\item [Type $E_6$] $r(\lieg)= 16$.
\item [Type $E_7$] $r(\lieg)= 27$.
\item [Type $E_8$] $r(\lieg)=57 $.
\item [Type $F_4$] $r(\lieg)= 15$.
\item [Type $G_2$] $r(\lieg)= 5$.
\end{description}
\end{example}
We note that for all root systems above, the minimal resonant codimension $r(\lieg)$ corresponds to the codimension of the maximal parabolic subalgebra $\lieq_{\Pi\sm\{\alpha_1\}}$ where the simple roots are as enumerated as in the Dynkin diagrams in Table \ref{table:bull}.

For the remainder of this subsection,  we show  that certain subgroups of $G$ with resonant codimension at most
$r(G)$ are parabolic.
With $\lieg$   the Lie algebra of $G$, let $\Sigma= \Sigma(\lieg)$ be the restricted root system of $\lieg$, and  let $$\lieg = \liec \oplus \liea \oplus \bigoplus_{\beta\in \Sigma} \lieg^\beta$$ be the restricted root space decomposition (relative to the choice of  Cartan involution $\theta$.)  Note that $\lieg^\beta$ is not a Lie subalgebra  if $2\beta$ is a root; in this case let $\langle \lieg^{\beta}\rangle$ denote the Lie-subalgebra generated by $\lieg^\beta$.

\begin{lemma}\label{lem:LieFacts1}
For any root $\beta\in \Sigma$, the subalgebra  $\liec$ acts irreducibly under the adjoint action on the root space $\lieg^{\beta}$.
\end{lemma}
As a corollary, let  $\lieh\subset \lieg$ be a Lie subalgebra with $\liec  \subset \lieh$.  Then for every $\beta \in \Sigma $ with  $\lieh \cap \lieg^\beta \neq 0$, we have $$\langle \lieg^{\beta}\rangle\subset \lieh.$$

\noindent A proof of Lemma \ref{lem:LieFacts1} using the complexification of $\lieg$ appears in  \cite[Lemma 5.3]{MR3748688}.  We give an alternative shorter proof of this fact using representation theory.
\begin{proof}[Proof of Lemma \ref{lem:LieFacts1}]
Let $V\subset  \lieg^{\beta}$ be a non-trivial, $\liec$-invariant subspace. Let $$\lieh = V\oplus \liec \oplus \liea.$$  Since the adjoint action of $\liea$  on $\lieg^\beta$ is by scalar multiplication and since $\liea$ centralizes $\liec$, it follows that $\lieh$ is a subalgebra.


 Fix a non-zero  $X\in V$.  By \cite[Lemma 7.73b]{MR1920389} applied to $X$ instead of $\theta(X)$, we have that
$$(\ad X)\colon (\ad X)(\lieg^{-\beta}  )\to  \lieg^{\beta} $$ is a bijection and since
$(\ad X)(\lieg^{-\beta})\subset  \lieg^0 = \liec\oplus \liea\subset \lieh$ it follows that $ \lieg^{\beta}\subset \lieh$.  
\end{proof}

\begin{proposition}\label{prop:LieFacts2}Let $\lieh\subset \lieg$ be a Lie subalgebra with $\liec \oplus \liea \subset \lieh$.
  If the cardinality of the set $\{[\beta]\in \hat \Sigma (\lieg) \colon \lieg^{[\beta]}\not \subset \lieh\}$ is at most   $r(\lieg)$, then $\lieh$ is parabolic.
\end{proposition}
Before we give the proof of Proposition \ref{prop:LieFacts2}, we need the following lemma whose proof requires case-by-case analysis.
In the analysis in the following lemma, we fix an inner product on $\liea^*$  which is preserved by the Weyl group and an orthonormal basis $\{e_1, e_2,\dots,\}$ for $\liea^*$ relative to which we may express all roots in a standard presentation such as in \cite[Appendix C] {MR1920389}.  Relative to the inner product, we may measure the lengths of roots.  All roots of the same length are in the same  orbit of the Weyl group.  If $\lieg$ is simple and if $\Sigma(\lieg)$ is of type $A_\ell, D_\ell, E_6, E_7,$ or $E_8$ all roots have the same length; if $\Sigma(\lieg)$ is of type $B_\ell, C_\ell, G_2,$ or $F_4$ there are two distinct lengths of roots and if $\Sigma(\lieg)$ is of type $(BC)_\ell$ there are three distinct lengths of roots.

\begin{lemma}\label{lem:disaster}
Let $\lieg$ be a simple Lie algebra and let  $\lieh\subset \lieg$ be a Lie subalgebra satisfying the hypotheses of Proposition \ref{prop:LieFacts2}.  Then either $\lieh= \lieg$ or there exists a long root $\beta_0$ such that $\lieg^{\beta_0}\cap \lieh = \{0\}$.
\end{lemma}
\begin{proof}
First note from Lemma \ref{lem:LieFacts1} that  $\lieh$ is saturated by full root spaces; that is, $\lieg^\beta \cap \lieh = \{0\}$ or $\lieg^\beta\subset \lieh$ for all roots $\beta\in \Sigma(\lieg)$.  If $\Sigma(\lieg)$ is of type $A_\ell, D_\ell, E_6, E_7,$ or $E_8$ then all roots are of the same length.  We  argue the lemma from case-by-case for the remaining possible abstract root systems $\Sigma(\lieg)$.

\subsection*{$\Sigma (\lieg)$ is of type $B_\ell$:} Relative to a choice of orthonormal basis $\{e_1, \dots , e_\ell\}$ the roots are $\{\pm e_i \pm e_j: 1\le i<j\le \ell\}\cup \{ \pm e_i:1\le i\le \ell\}$; the long roots are $\{\pm e_i \pm e_j\}$.  Suppose $\lieg^\beta\subset \lieh$ for all long roots $\beta$.  Since $r(\lieg)= 2\ell-1$, by hypotheses and assumption there exists $1\le i_0\le \ell$ and a short root $\beta'\in \{\pm e_{i_0}\}$ with $\lieg^{\beta'}\subset \lieh$.  For each $1\le i\le \ell$ we have  $\lieg^{\pm e_i-e_{i_0}}, \lieg^{\pm e_i+ e_{i_0}} \subset \lieh$.  Bracketing with $\lieg^{\beta'}$, we have  $\lieg^{\pm e_i}\subset \lieh$ for every $1\le i\le \ell$.   It follows that $\lieh= \lieg$.

\subsection*{$\Sigma (\lieg)$ is of type $C_\ell$:} The roots are  $\{\pm e_i \pm e_j: 1\le  i<j\le \ell\}\cup \{ \pm 2 e_i\}$; the long roots are $\{\pm 2 e_i\}$.  We induct on $\ell$.  In the case $\ell=2$, the conclusion follows from the above since $C_2$ and $B_2$ are isomorphic.  Suppose $\lieg^\beta\subset \lieh$ for all long roots $\beta\in \{\pm 2 e_i\}$.  For the sake of contradiction, assume $\lieh\neq \lieg$; then there is a short root $\beta'= \pm e_i \pm e_j$ with $\lieg^{\beta'}\cap \lieh = \{0\}$.  Acting by the Weyl group, we may assume $\beta'=e_1-e_2 = \alpha_1$ is the left-most root in the Dynkin diagram with respect to some base $\Pi$.  Let $\lieg'$ be the Lie subalgebra of $\lieg$ generated by the root spaces associated to roots $\pm \alpha_2, \dots, \pm \alpha_{\ell}$.  Then $\Sigma(\lieg')$ is of type $C_{\ell-1}$.

Since $\lieg^{e_1-e_2}\cap \lieh =\{0\}$ and $\lieg^\beta\subset \lieh$ for $\beta \in\{\pm2 e_1, \pm 2e_2\}$, we conclude that
$\lieg^{\beta'}\cap \lieh =\{0\}$ for the 4 roots $\beta'\in \{\pm e_1 \pm e_2\}.$
Let $\lieh'= \lieh\cap \lieg'$.  Then  the cardinality of the set  $\{[\beta]\in \hat \Sigma (\lieg') \colon \lieg^{[\beta]}\not \subset \lieh'\}$ is at most $r(\lieg)-4= 2\ell-1-4=2(\ell-1)-3<r(\lieg')$.  In particular, $\lieh'\subset \lieg'$ satisfies the hypotheses of Proposition \ref{prop:LieFacts2} and, since $\lieg'$ contains all root spaces associated to its long roots, by the inductive hypotheses we conclude that $\lieh'=\lieg'$.

Finally, there are $4\ell-4$ roots of the form $\pm e_1 \pm e_j$, $j\ge 2$.  As we assume $\lieh$ contains all root spaces associated to long roots and since $r(\lieg) =2\ell-1<4\ell-4$, there exists $2\le i_0\le \ell$ such that $\lieg^{\beta'}\subset \lieh$ for some, and hence all, roots  $\beta' \in \{\pm e_1 \pm e_{i_0}\}$.  Since $\lieg'\subset\lieg$ and since $\pm e_1 \pm e_j = (\pm e_1 - e_{i_0}) + (e_{i_0} \pm e_j )$,   we conclude that
$\lieg^{\beta'}\subset \lieh$ for all  $\beta' \in \{\pm e_{1} \pm e_{j}: 2\le j\le \ell\}$.  It follows that $\lieh= \lieg$, contradicting the assumption $\lieh\neq \lieg$ above.

\subsection*{$\Sigma (\lieg)$ is of type $(BC)_\ell$:}  The roots are  $\{\pm e_i \pm e_j\}\cup \{ \pm e_i\}\cup \{ \pm 2 e_i\}$; the long roots are $\{\pm 2 e_i\}$.
Suppose $\lieg^\beta\subset \lieh$ for all long roots $\beta\in \{\pm 2 e_i\}$.  From the previous analysis, $\lieh$ contains the subalgebra (with root system of type $C_\ell$) containing all root spaces $\lieg^\beta$ associated to roots of the form $\beta = \{\pm e_i \pm e_j\}\cup \{\pm 2e_i\}$.
From the analysis when $\Sigma (\lieg)$ is of type $B_\ell$, it follows that $\lieg^{\beta'}\subset \lieh$ for every root $\beta'\in \{\pm e_j \}$  and thus $\lieh=\lieg$.

\subsection*{$\Sigma (\lieg)$ is of type $G_2$:}
The roots are $\{e_i-e_j: 1\le i,j\le 3: i\neq j\} \cup
\{\pm (2e_1-e_2-e_3),\pm (2e_2-e_1-e_3), \pm (2e_3-e_1-e_2)\}$; the long roots are  $\{\pm (2e_1-e_2-e_3),\pm (2e_2-e_1-e_3), \pm (2e_3-e_1-e_2)\}$.
Suppose $\lieg^\beta\subset \lieh$ for all long roots $\beta$.  Since $r(\lieg) = 5$ there is at least one short root $\beta'$ with $\lieg^{\beta'}\subset \lieh$; acting by the Weyl group, we may assume $\beta'= e_1-e_2$.
Observe that  $(e_2-e_3) = 2e_2-e_1-e_3 + \beta'$,
$(e_1-e_3) = \beta' + (e_2 - e_3)$,  $(e_3-e_1) = -(2e_1-e_2-e_3) + \beta'$,
$(e_3-e_2)=(e_3-e_1) + \beta'$, and
$(e_2-e_1)=(e_2-e_3) + (e_3-e_1)$.  It follows that $\lieg^{\beta}\subset \lieh$ for all short roots $\beta$ whence $\lieh= \lieg$.

\subsection*{$\Sigma (\lieg)$ is of type $F_4$:}
There are 48 roots  $\{\pm e_i \pm e_j: 1\le i <j\le 4\} \cup \{\pm e_i: 1\le i\le 4\} \cup \{\frac 1 2 (\pm e_1 \pm e_2\pm e_3\pm e_4)\}$;  the long roots are the 24 roots $\{\pm e_i \pm e_j: 1\le i <j\le 4\}$.  Suppose $\lieg^\beta\subset \lieh$ for all long roots $\beta$.  We have  $r(\lieg) = 15$ so there is at least one short root $\beta'$  with $\lieg^{\beta'}\subset \lieh$; acting by the Weyl group, we may assume $\beta'= \frac 1 2 ( e_1  + e_2+ e_3 +e_4)$.
Taking brackets of $\lieg^{\beta'}$ with $\lieg^\beta$ for all long roots $\beta$, we have
$\lieg^{\beta''}\subset \lieh$  for the 8 roots $\beta''\in  \{\frac 1 2 (\pm e_1 \pm e_2\pm e_3\pm e_4)\}$ with an even number of positive terms.

We also claim $\lieg^{\wtd \beta}\subset \lieh$ for at least one $\wtd \beta = \frac 1 2 (\pm e_1 \pm e_2\pm e_3\pm e_4)$ with an odd number of positive terms.  Indeed  there are 8 roots of the form  $\frac 1 2 (\pm e_1 \pm e_2\pm e_3\pm e_4)$ with an odd number of positive terms and 8 roots of the form $\pm e_i$. Since $r(\lieg) = 15$, one of these 16 roots corresponds to a root space in $\lieh$; the former case gives such a $\wtd \beta$  and the latter case gives such $\wtd \beta$ after bracketing with some root space   associated to a root $  \frac 1 2 (\pm e_1 \pm e_2\pm e_3\pm e_4)$ with an even number of positive terms.  Taking brackets of $\lieg^{\wtd \beta}$ with $\lieg^\beta$ for all long roots $\beta$, we have $\lieg^{\bar \beta}\subset \lieh$  for the 8 roots $\bar \beta \in  \{\frac 1 2 (\pm e_1 \pm e_2\pm e_3\pm e_4)\}$ with an odd number of positive terms.

Finally,  we have $ e_1 = \frac 1 2 ( e_1 + e_2+ e_3+  e_4) + \frac 1 2 ( e_1 - e_2- e_3- e_4)$ and
$- e_1 = \frac 1 2 (- e_1 + e_2+ e_3+  e_4) + \frac 1 2 (- e_1 - e_2- e_3- e_4)$.
Moreover, for $2\le i\le 4,$ we have $\pm e_i = (-e_1 \pm  e_i ) + e_1$.  It follows that  $\lieg^{\beta'''}\subset \lieh$ for  the 8 roots $\beta'''\in \{\pm  e_i\}$.  Combined with  the above analysis, it follows that $\lieh = \lieg$.
\end{proof}

\begin{proof}[Proof of Proposition \ref{prop:LieFacts2}]  
First recall from Remark \ref{rem:ooppllkk} that  $r(\lieg)$ is the minimal value of $r(\lieg')$ as $\lieg'$ varies over simple non-compact ideals in $\lieg$.  In particular, if the conclusion holds for all simple Lie algebras $\lieg$, then it automatically holds for all semisimple Lie algebras.  Thus we may assume $\lieg$ is simple for the remainder.

We may assume $\lieh\neq \lieg$.  Let $\lieh'\subset \lieh$ be the Lie subalgebra generated by $\liec$, $\liea$, and all coarse root spaces $\lieg^{[\beta]}$ where $\lieg^{[\beta]}\subset \lieh$.
It follows from Lemma \ref{lem:disaster} that there exists a long root $\beta_0$ with $\lieg^{\beta_0}\cap \lieh = \{0\}$.  Acting by the Weyl group, the root $-\beta_0$ is the highest root for some choice of base $\Pi$.  If $\Pi = \{\alpha_1, \dots, \alpha_\ell\}$ are the simple positive roots for this base, select $X\in \liea$ such that $\alpha_j(X)> 0$  for all $1\le i\le \ell$.
Then $\beta_0(X)$ is the minimal value of $\beta(X)$ as $\beta$ varies over all $\beta\in \Sigma(\lieg)$; moreover for all $\beta\in \Sigma(\lieg)\sm \{\beta_0\}$ we have $\beta_0(X)< \beta(X).$

Let $\lien$ be the Lie subalgebra generated by all positive roots relative to $\Pi$.  Let $\liel$ be the Lie subalgebra generated by $\lieh'$ and $\lien$; observe that $\liel$ is parabolic.   Since the minimal value of $\beta(X)$ is $\beta_0(X)$,   it follows that $\lieg^{\beta_0}$ is not a subspace of $\liel$; in particular $\liel\neq \lieg$.
Since $\lieh'\subset \liel\neq \lieg$, it follows from the definition of $r(\lieg)$ that $\lieh' = \liel$.
Since $\liel=\lieh'\subset \lieh$, the conclusion follows.
\end{proof}



For the reader less familiar with finite dimensional representation theory and root systems
we also include a geometric proof of a weaker assertion that suffices for all our proofs in the
case of $\R$-split groups.  It is also possible to give a proof of Proposition \ref{prop:LieFacts2} above
using Lemma \ref{lem:dave} below and Lemma \ref{lem:LieFacts1}.

\begin{lemma}\label{lem:dave}
Let $\lieh\subset \lieg$  be a subalgebra whose codimension is at most the minimal codimension of all proper parabolic subalgebras of $\lieg$.  Then $\lieh$ is parabolic.
\end{lemma}
\begin{proof}
We may assume that $\dim \lieh$ is maximal among all proper subalgebras of~$\lieg$. Let $H$ be the connected Lie subgroup of~$G$ whose Lie algebra is~$\lieh$. As the statement concerns  Lie algebras, we may replace $G$ with its adjoint group and assume that $G$ is a (real) linear algebraic group. So we may let $\overline{H}$ be the Zariski closure of~$H$ in~$G$. Since $H$ is connected, we know that $H$ is not Zariski dense, so $\dim \overline{H} < \dim G$. Then the maximality of $\dim H$ implies that $H$ is the identity component of~$\overline{H}$, and therefore has finite index in~$\overline{H}$ so $\mathfrak{h}$ is also the Lie algebra of the real algebraic group~$\overline{H}$.

By Chevalley's Lemma \cite[Proposition 3.1.4]{MR0776417}, there is a finite-dimensional representation $\rho \colon G \to \Gl(n,\R)$, such that $\overline{H}$ is the stabilizer of a point~$x$ in the corresponding projective space $\R P^{n-1}$. Since finite-dimensional representations of~$G$ are completely reducible and $G$ has no nontrivial $1$-dimensional representations, we may assume without loss of generality that $G$ has no fixed points in $\R P^{n-1}$.

Since this is an action of an algebraic group on a variety, we know that the closure of the $G$-orbit of~$x$ consists of the union of $Gx$ with orbits of strictly smaller dimension. However, the maximality of $\dim H$ and the absence of fixed points implies that there are no $G$-orbits of smaller dimension. So $Gx$ must be a closed subset of $\R P^{n-1}$, and is therefore compact. This means $G/\overline{H}$ is compact. 

If $H$ and $\overline{H}$ are reductive, then they are unimodular and $G/\overline{H}$ admits a finite invariant measure.  By the Borel density theorem \cite[Theorem 3.2.5]{MR0776417} this implies $G=\overline{H}$, a contradiction.  If $\overline{H}$ is not reductive then the unipotent radical, $U$, of~$\overline{H}$ is nontrivial. A result of Borel and Tits \cite[Proposition 3.1]{MR294349} states that $U$ is contained in the unipotent radical of a parabolic subgroup~$P$ that contains the normalizer of $U$. Since $\overline{H}$ is contained in this normalizer, it must be contained in~$P$.  Moreover $P$ is proper because its unipotent radical contains $U$, and is therefore nontrivial. Then the maximality of $\dim H$ implies that $H$ is the identity component of~$P$, so $\lieh$ is the Lie algebra of~$P$. This is also a consequence of the more detailed result \cite[Thm.~1.2]{MR1079114}.
\end{proof}

\section{Suspension action and Proof of Theorem \ref{prop:USEGOD}}
 We begin  by introducing the suspension action with which we work for the remainder of the proof of  Theorem \ref{prop:USEGOD}.  We then give some general background on  Lyapunov exponents and state the two key propositions used in the proof of Theorem \ref{prop:USEGOD}.

\subsection{Suspension space}
\label{subsection:suspension}
\label{sec:susp}
Recall we fix $G$ to be a semisimple Lie group with real-rank at least $2$.  Let $\Gamma\subset G$ be a cocompact lattice and   let $\alpha \colon \Gamma\to \diff^{1+\beta}(M)$ be an action for $\beta>0$.

We construct an auxiliary   space  on which the  action $\alpha$ of $\Gamma$ on $M$ embeds as a Poincar\' e section for an associated $G$-action.
On the product $G\times M$ consider the right $\Gamma$-action
	$$ (g,x)\cdot \gamma= (g\gamma , \alpha(\gamma\inv)(x))$$ and the left $G$-action $$a\cdot (g,x) = (ag, x).$$
Define  the quotient manifold $M^\alpha:= G\times M/\Gamma $.  As the  $G$-action on $G\times M$ commutes with the $\Gamma$-action, we have an induced left $G$-action  on $M^\alpha$.  For $g\in G$ and $x\in M^\alpha$ we denote this action by $g\cdot x$ and denote the derivative of the diffeomorphism $x\mapsto g\cdot x$ by $Dg$.
We write $\pi\colon M^\alpha\to G/\Gamma$ for the natural projection map.  Note that $M^\alpha$ has the structure of a fiber-bundle over $G/\Gamma$ induced by the map $\pi$ with fibers diffeomorphic to $M$. Note that the $G$-action preserves the fibers of $M^\alpha$.
  As the action of $\alpha$ is by $C^{2}$ diffeomorphisms, $M^\alpha$ is naturally a  $C^{2}$ manifold.  Equip $M^{\alpha}$ with a $C^\infty$ structure compatible with the $C^{2}$ structure; the existence of this compatible structure is guaranteed by a classical theorem of Whitney \cite[Theorem 1]{MR1503303}.
Choose a right-$\Gamma$-invariant Riemannian metric on $G\times M$ whose restriction to any $G \times \{m\}$ is right-$G$-invariant. This exists because the $\Gamma$-action on $G \times M$ is proper.  This metric defines a Riemannian metric on $M^{\alpha}$ whose restriction to the tangent space to any $G$-orbit pushes forward to a metric on $G/\Gamma$ defined by a right-$G$-invariant metric on $G$.

\subsection{Lyapunov exponents and Oseledec's theorem}
\label{subsection:lyapunov}
Let $X$ be a compact metric space equipped with a continuous (left) $G$-action.  
A measurable function $\calA \colon G\times X\to \Gl(d,\R)$ defines a linear cocycle if
$$\calA (g',g\cdot x)\calA(g,x) = \calA (g'g, x).$$
Then $\calA$ defines an action by vector bundle automorphisms on the trivial bundle $X\times \R^d$ which projects to the $G$-action on $X$.
%

More generally, let $\calE\to X$ be a continuous, finite dimensional, normed vector bundle.  A measurable linear cocycle over the $G$-action on $X$ is a measurable action $$\calA\colon G\times \calE\to \calE$$ by vector-bundle automorphisms that projects to the $G$-action on $X$.
We write $\calE_x$ for the fiber of $\calE$ over $x$ and  $$\calA(g,x)\colon \calE_x\to\calE_{g\cdot x} $$ for the linear map between Banach spaces $\calE_x$ and $\calE_{g\cdot x}$.

Below, we will always assume our cocycle $\calA \colon G\times \calE\to \calE$ is bounded:  for every compact $K\subset G$
$$\sup_{(g,x)\in K\times X} \|\calA(g,x)\|$$ is bounded.
Moreover, we typically assume the  action $\calA \colon G\times \calE\to \calE$ is continuous which automatically implies boundedness.   If one cares only about measurable cocycles, one may assume the bundle $\calE$ is trivial.

%

Given $s\in G$ and an $s$-invariant Borel probability measure $\mu$ on $X$ we define the \emph{average top} (or \emph{leading})  \emph{Lyapunov exponent of $\calA$} to be
\begin{equation} \label{eq:toplyap} \lambda_+(s,\mu, \calA) := \inf _{n} \frac 1 n  \int \log \|\calA(s^n, x)\| \ d \mu (x).\end{equation}
We also note that for an $s$-invariant measure $\mu$,  the sequence  $\frac 1 n  \int \log \|\calA(s^n, x)\| \ d \mu (x)$ is subadditive whence the infimum in \eqref{eq:toplyap} may be replaced by a limit. By the Kingman subadditive ergodic theorem (see \cite[Theorem 3.3]{viana2014lectures}) if $\mu$ is ergodic, the functions $\frac 1 n    \log \|\calA(s^n, x)\|$ converges $\mu$-\ae to $\lambda_+(s,\mu, \calA)$ as $n \rightarrow \infty$.

We have the following elementary fact.
\begin{claim}\label{claim:jjjiiilllkkk}
If the restriction of the cocycle to $\calA\colon G\times \calE\to \calE$ to $s\in G$ is continuous then the map $$\mu \mapsto \lambda_+(s,\mu, \calA)$$ is upper-semicontinuous on the set of all $s$-invariant Borel probability measures equipped with the weak-$*$ topology.
\end{claim}

We  recall the following standard fact which is crucial in our  later  averaging arguments.  
Given an amenable subgroup $H\subset G$, a bounded measurable set $F\subset H$ of positive Haar measure, and a probability measure $\mu$ on $X$ denote by $F\ast \mu$   the probability measure defined as follows: for a Borel $B\subset X$ let   $$(F\ast \mu) (B)  = \dfrac 1 {|F|}\int_F \mu (s\inv\cdot B)\ d s$$ where $|F|$ is the volume of the set $F$ induced by the (left-)Haar measure on $H$.
For $x\in X$, we write  $$\nu_x^{F} = F\ast \delta_x .$$


\begin{lemma}\label{lemma:averagingiscool}
Let $\calA\colon G\times \calE\to \calE$ 
be a bounded, continuous linear cocycle.
Let $s\in G$ and let $\mu$ be an $s$-invariant, Borel probability measure on $X$.
Let $H\subset G$ be an amenable subgroup contained in the centralizer of $s$ in $G$.  Let $F_{m}$ be a \Folner sequence of precompact sets in $H$ and let $\mu'$ be a Borel probability measure that is a weak-$*$ subsequential limit of the sequence of measures $\{F_{m} \ast \mu\}$. Then
\begin{enumlemma}
\item \label{averagingiscool1}$\mu'$ is $s$-invariant and $H$-invariant;
\item \label{averagingiscool2}$\lambda_+(s,\mu', \calA)  \ge \lambda_+(s,\mu, \calA) $.
\end{enumlemma}
\end{lemma}

\begin{proof}
 \ref{averagingiscool1} follows as each $\{F_m \ast \mu\}$ is $s$-invariant and $s$-invariance is closed under weak-$*$ convergence.

  For \ref{averagingiscool2},  first note that as $\calA$  is assumed bounded, it follows from the cocycle relation that  $\lambda_+(s,F_m \ast \mu, \calA) =  \lambda_+(s,\mu, \calA)$ for every $m$.
Indeed, for any $t \in H$ let  $C_t = \sup_{x \in X} \log \|\calA(t^{\pm 1}, x)\| $ and let $C_m = \sup _{t\in F_m } C_t$.  For  $x \in M$ and   $t \in F_m$ the cocycle property and subadditivity of norms yields
\begin{align*}
\log \|\calA(s^n, tx)\| &\leq  C_t + \log \|\calA(s^nt, x)\|  \\  
&=  C_t + \log \|\calA(ts^n, x)\|  \\
&\leq 2C_t + \log \|\calA(s^n, x)\|\\
&\leq 2C_m+ \log \|\calA(s^n, x)\|.
\end{align*}
Similarly we can prove that  $\log \|\calA(s^n, tx)\| \geq  \log \|\calA(s^n, x)\|  - 2C_m$.\\
%
Thus,
\begin{align*}
\int \log \|\calA(s^n, x) \| \ d (F_m\ast \mu )(x) &= \dfrac 1 {|F_m|}\int_{F_m} \int \log \|\calA(s^n, x)\| \ dt*\mu (x) \\
&= \dfrac 1 {|F_{m}|}\int_{F_{m}} \int \log \|\calA(s^n, tx)\| \ dt d\mu (x) \\
&\leq \dfrac 1 {|F_{m}|}\int_{F_{m}}  \bigg( \ \ 2C_m + \int \log \|\calA(s^n, x)\|  \ d \mu (x) \bigg) \ d t \\
&\leq 2C_m + \int \log \|\calA(s^n, x)\| d \mu{(x)}
\end{align*}
Dividing by $n$ yields $\lambda_+(s,F_{m} \ast \mu, \calA) \leq  \lambda_+(s,\mu, \calA)$. The reverse  inequality is similar.   
Conclusion \ref{averagingiscool2}      follows  from the upper semicontinuity in Claim \ref{claim:jjjiiilllkkk}.
\end{proof}

Consider $A\subset G$   any abelian subgroup isomorphic to $\R^k$.  Equip $A\cong \R^k$ with any norm $|\cdot|$.  Consider   an $A$-invariant, $A$-ergodic probability  measure $\mu$ on $X$.  For a bounded measurable linear cocycle $\calA\colon A\times X\to \Gl(d,\R)$ we have the following consequence of the higher-rank Oseledec's multiplicative ergodic theorem  (c.f.\ \cite[Theorem 2.4]{AWB-GLY-P1}).
\begin{proposition}
\label{thm:higherrankMET}
There are
	\begin{enumerate}
	\item an $\alpha$-invariant subset $\Lambda_0\subset X$ with $\mu(\Lambda_0)=1$;
  \item 
   linear functionals $\lambda_i\colon \R^k\to \R$ for $1\le i\le p$;  
	\item   and splittings   $\R^d= \bigoplus _{i=1}^p E_{\lambda_i}(x)$ 
	into families of mutually transverse,  $\mu$-measurable  subbundles $E_{\lambda_i}(x)\subset \R^d$ defined  for $x\in \Lambda_0$

	\end{enumerate}
such that
\begin{enumlemma}	
	\item $\calA (s, x) E_{\lambda_i}(x)= E_{\lambda_i}(s\cdot x)$ and
	\item \label{lemma:partb} $\displaystyle \lim_{|s|\to \infty} \frac { \log |  \calA (s,x) (v)| - \lambda_i(s)}{|s|}=0$
\end{enumlemma}	
	for all $x\in \Lambda_0$ and all $ v\in  E_{\lambda_i}(x)\sm \{0\}$.  
 \end{proposition}
Note that \ref{lemma:partb} implies for $v\in E_{\lambda_i}(x)$ the weaker result that for $s\in A$,

$$\lim_{k\to\pm \infty} \tfrac {1} k \log |  \calA (s^k,x) (v)| =  \lambda_i(s)$$

We also remark that if $\mu$ is an $A$-invariant, $A$-ergodic measure  then for any $s\in A$ the average top Lyapunov exponent is given as

\begin{equation}\label{eqn:lyapn}
\lambda_+(s,\mu, \calA) = \max _i \lambda_i(s).
\end{equation}

In the case that $\mu$ is $A$-invariant but not $A$-ergodic, Proposition \ref{thm:higherrankMET} holds on each $A$-ergodic component of $\mu$. Even more is true, the number of Lyapunov exponents is determined by an integer valued measurable function $1 \leq p(x) \leq k$ constant on ergodic components and all the data arising from Proposition \ref{thm:higherrankMET}, including the linear functionals and the subspaces, varies measurably in $X$; see \cite[Section 3.6.1]{MR2348606}.  In this case we have the following construction which will be used later to avoid passing to ergodic components.

\begin{lemma}\label{lemmatoskipergodicity}
If $\mu$ is an $A$-invariant measure  on $X$ then for any $s'\in A$, there is a linear functional $\lambda_{+,s',\mu}\colon A\to \R$ so that
\begin{enumerate}
\item  $\lambda_{+,s',\mu}(ts') =  \lambda_+(ts',\mu, \calA)$ for any $t\ge 0 $;
\item $\lambda_+(s,\mu, \calA)\ge \lambda_{+,s',\mu}(s)$ for all $s\in A$.
\end{enumerate}
\end{lemma}

\begin{proof}
We first pass to  an ergodic decomposition of the $A$-action on $(X,\mu)$.  See for instance \cite[Theorem 2.19]{MR2114789}.  This gives a Borel map $\zeta \colon X \rightarrow \Omega$ where
$\Omega$ is the space of ergodic components of the $A$ action and a Borel map $\xi\colon \Omega \rightarrow \Prob(X)$ where the target is the space of probability measures on $X$,
such that
$$\mu= \int_{\Omega} \xi(\omega) d\zeta_*\mu.$$
See for example \cite{MR2114789} for more details. Since the function $p$ mentioned in the paragraph preceding this lemma is constant on ergodic components, we
can view it as a function on $\Omega$.

By the dominated convergence theorem, one checks that 
$$\lambda_+{(s,\mu, \calA)} = \int_{\Omega} \lambda_+{(s, \xi(\omega), \calA)} \ d \zeta_*\mu.$$
 From this and  \eqref{eqn:lyapn}, one verifies that  $\lambda_+(ts,\mu, \calA)=t\lambda_+(s,\mu, \calA)$ for any positive real number $t$.  This can also be proven directly from properties of $\lambda_+(s,\mu, \calA)$.  Once we construct a linear functional $\lambda_{+,s',\mu}$, this immediately implies the first claim of the lemma.

Let $\calL_{\omega}= \{\lambda_{i,\omega}, 1\le i \le p(\omega)\}$ be the collection of   Lyapunov exponents for the cocycle $\calA$ and the measure $\xi(\omega)$.
Let $\chi\colon \omega \mapsto \lambda_{+,s',\omega} \in \calL_{\omega}$ be a measurable, $A$-invariant assignment  satisfying
$$\lambda_{+,s',\omega}(s') = \max_i \lambda_{i,\omega}(s').$$  We briefly defer justifying the existence of $\chi$.
Take $\lambda_{+,s',\mu}\colon A\to \R$ to be
	$$\lambda_{+,s',\mu}(s) = \int   \lambda_{+,s',\omega}(s) \ d \zeta_*\mu.$$
The integral is defined since $\lambda_{+,s',\omega}(s)$ is bounded and one verifies that  $\lambda_{+,s'}(s) $ satisfies the properties of the Lemma.

To justify the existence of the measurable map $\chi$ one can use the measurable selection theorem (see e.g. \cite[Theorem 2.3]{MR2114789}), but one can also give
a simpler argument. We construct this map from a map $\chi_X$ from $X$ that factors through $\Omega$.  Since we can partition $X$ into finitely many disjoint measurable
subsets where $p(x)$ is constant, we assume $p=p(x)$ is constant.
Let $X_p$ be the union of $p$ disjoint copies of $X$ and let $\calE^*$ be the dual bundle to $\calE$.  There is a measurable map $\chi_p\colon X_p \rightarrow \calE^*$
sending $x$ to the set $\calL_{\zeta(x)}$.  Choosing $s' \in A$, we can define a subset of $X^{\max}$ and a restriction $\chi_{\max}\colon X^{\max} \rightarrow \calE^*$
where $X^{\max}$ consists of those linear functional in $\calL_{\zeta(x)}$ such that $\lambda_{i,\omega}(s')=\max_i \lambda_{i,\omega}(s')$.  We can
partition $X$ into finitely many measurable sets $X_i$ where $X^{\max}$ is $i$ disjoint copies of $X$ for $i \leq p$.  On each $X_i$ we can make a choice of one copy
of $X_i$ which choses the linear functional which is the image of $\chi$ on $X_i$.  This assignment is clearly measurable on $X_i$.
\end{proof}

\subsection{Subexponential growth of fiberwise derivatives}
We return to the setting introduced in Section \ref{sec:susp}.
With $\pi\colon M^\alpha\to G/\Gamma$  the projection,  let $F = \ker (D\pi)$ denote  the fiberwise tangent bundle of $M^\alpha$.

We say the induced action of $G$ on $M^\alpha$ has \emph{uniform subexponential growth of fiberwise derivatives} if for all $\epsilon>0$ there is a $C$ such that $$\|\restrict{Dg} {F}\|\le Ce^{\epsilon d(e, g)}$$
where  $\|\restrict{Dg} {F}\|=\sup_{x\in M^{\alpha}}\|\restrict{Dg(x)} {F(x)}\|$. As $\Gamma$ is cocompact,  there is a clear relation between the  growth of the fiberwise derivatives for the $G$-action and the growth of derivatives of the $\Gamma$-action.

\begin{claim}\label{claim:gammavG}
The action $\alpha$ of $\Gamma$ on $M$ has uniform subexponential growth of derivatives if and only if the induced action of $G$ on $M^\alpha$ has uniform subexponential growth of  fiberwise derivatives.
\end{claim}

\subsection{Proof of Theorem \ref{prop:USEGOD}}
We let $\calA$ denote the fiberwise derivative cocycle for the action of $G$ on $M^\alpha$; that is $\calA(g, x) = \restrict {D_x g} { F}$.  Let $A= \exp \liea \subset G$ be a maximal split Cartan subgroup.  Given $s\in A$ and an $s$-invariant Borel probability measure $\mu$ we
write $$\lambda^F_+(s,\mu) := \lambda_+(s,\mu, \calA)  =  \inf _{n\to \infty} \frac 1 n  \int \log \|\restrict{D_x(s^n)}{F} \| \ d \mu (x)$$
for the \emph{average top fiberwise Lyapunov exponent} of $s$ with respect to $\mu$.

The proof of Theorem \ref{prop:USEGOD} is by contradiction.  Assuming  Theorem \ref{prop:USEGOD}  fails, from Claim \ref{claim:gammavG} we first establish the following.

\begin{proposition}\label{prop:lyapons}
Suppose the induced action of $G$ on $M^\alpha$ fails to have uniform subexponential growth of  fiberwise derivatives.
Then there is an $s\in A$ and an $A$-invariant Borel probability measure $\mu$ with  $\lambda^F_+(s,\mu) >0$.  
\end{proposition}

As discussed above,  Theorem \ref{prop:USEGOD} holds trivially in the case where $G$ has rank-1 factors.   To complete the proof of Theorem \ref{prop:USEGOD}  we may thus assume
that all non-compact, almost-simple factors of $G$ are of higher-rank.   The proof of the following proposition contains the major technical innovations in this paper.

\begin{proposition}\label{prop:Ginvariantmeasure}
Let $G$ be a connected, semisimple Lie group with finite center, all of whose non-compact, almost-simple factors are of real-rank at least 2.   Let  $\Gamma\subset G$ be a cocompact lattice and let $\alpha\colon \Gamma\to \diff^{1+\beta}(M)$ be an action.  Suppose that either
\begin{enumerate}
\item $\dim(M)<r(G)$, or
\item $\dim(M)= r(G)$ and $\alpha$ preserves a smooth volume
\end{enumerate}
and that there is an $s\in A$ and an $A$-invariant Borel probability measure $\mu$ on $M^\alpha$ with $\lambda_+^F(s,\mu)>0$.
Then there is a $G$-invariant measure $\mu'$ and $s'\in A$ with $\lambda_+^F(s',\mu')>0$.
 \end{proposition}

From Proposition \ref{prop:Ginvariantmeasure} we immediately obtain a contradiction with Zimmer's cocycle superrigidity theorem and the fact that there are no non-trivial linear representations of $G$ into $\Gl(r(G), \R)$ \cite{MR900826}. Theorem \ref{prop:USEGOD} follows immediately from Propositions \ref{prop:lyapons},  \ref{prop:Ginvariantmeasure} and Claim \ref{claim:gammavG}.


\section{Proof of Proposition \ref{prop:lyapons}}
\def\e{\epsilon}
To establish Proposition \ref{prop:lyapons}, suppose the induced action of $G$ on $M^\alpha$ fails to have uniform subexponential growth of  fiberwise derivatives.
Then there exist $\epsilon>0$, a sequence of elements $g_n$ in $G$ with $d(e, g_n)\to \infty$, a sequence of base points $x_n \in M^{\alpha}$, and a sequence of  unit vectors $v_n\in F(x_n):= T_{x_n}M^{\alpha} \cap F$ tangent to the fibers of $M_\alpha$ satisfying  $$\|D_{x_n} g_n (v_n)\| \geq e^{3\e d(e, g_n)}.$$

Let  $G = KAK$  be the Cartan decomposition of $G$  (c.f. \cite[Theorem 7.39]{MR1920389}).
For each $g_n$, write  $g_n= k_na_nk'_n$ where $k_n, k_n' \in K$ and $a_n \in A$. Note that $a_n\to \infty $ as $n\to \infty$.  As $K$ is a compact, the fiberwise derivative $\sup_{k\in K}\|\restrict{Dk} {F}\|$ is bounded above and thus  $$\|D_{x_n}a_n(v_n)\| \geq e^{2 \e d(a_n, e)}$$ for all sufficiently large $n$.

Recall that the Lie exponential $\exp\colon \lieg\to G$ restricts to a diffeomorphism from $\liea$ to $A$; moreover, $\exp \colon \liea \to A $ is an isometry.  Write $a_n = \exp (t_nu_n)$ where $u_n$ is a unit vector in $\liea$ and $t_n = d(a_n,e)$.   Given $t\in \R$, let $[t]$ denote the integer part of $t$.  Then for  sufficiently large $n$ we have
$$\|D_{x_n}\exp ([t_n]u_n)(v_n)\| \geq e^{ \e [t_n]}.$$
Passing to a subsequence, we  assume  $u_n$  converges to a unit vector $u\in \liea$.  The element $s = \text{exp}(u)\in A$ will be the element satisfying the conclusion of the proposition.

\def\P{\mathbb{P}}
Recall $F=\ker (D\pi)$ denotes the fiberwise tangent bundle of $M^\alpha$.  Let $UF$ denote associated the unit-sphere bundle; that is, the quotient  of $F$ under the equivalence relation of positive proportionality in each fiber $F(x)$ of $ F$.
We represent elements of $UF$ by pairs elements $(x,v)$ where $x \in M^{\alpha}$ and $v$ is a unit vector in the fiber  $F(x)$.
The derivative of the $G$-action on $M^\alpha$ induces a $G$-action on $F$ by fiber-bundle automorphisms;  the map intertwining fibers is denoted by $D_x g\colon F(x)\to F(gx)$.  The $G$-action of $F$ induces a $G$-action on $UF$; we denote the map intertwining fibers of $UF$ by $U D_x g\colon UF(x)\to UF(gx)$.

 For each   $n$, we define a  Borel probability measure $\nu_n$ on $UF$ as follows: Given a continuous $\phi\colon UF\to   \R$, let
$$\int \phi \ d \nu_n := \frac{1}{[t_n]}\sum_{m= 0}^{[t_n]-1} \phi \big(\exp(m u_n)\cdot (x), UD_x\exp(mu_n) ( v_n)\big) .$$
Given $g\in G$ and  a probability measure $\nu$ on  $UF$ consider  the expression $$\psi(g, \nu) = \int_{UF} \log \bigg( \frac{\|D_xg(v)\|_{gx}}{\|v\|_x} \bigg) \ d\nu(x,v).$$
From the definition of $\nu_n$   we have for   every $n$ that
\begin{equation}\label{facil}
  \psi({\exp}(  u_n), \nu_n) \geq   \e. 
\end{equation}

Consider any  weak-$*$ accumulation point    $\nu$   of the sequence of probability measures $\{\nu_n\}$ on $UF$.  
We have that $\nu$ is invariant under $s:= \exp(u)$.  
Indeed, let  $f\colon UF\to \R$ be a continuous function.  Then
\begin{align*}
\int_{UF} f \circ s-  f \ d \nu_n &= \int_{UF} f\circ \exp(u) - f \circ \exp(u_n) \ d\nu_n\\
&+ \int_{UF}  f \circ \exp(u_n) - f  \ d\nu_n
\end{align*}

The first integral converges to zero as the functions $f\circ \exp(u) - f \circ \exp(u_n)$ converges uniformly to zero in $n$.  The second integral clearly converges to zero  by compactness and the definition of $\nu_n$.
Taking $n\to \infty$, $$\int_{UTM^{\alpha}} f \circ s  \ d \nu = \int_{UTM^{\alpha}} f \ d \nu.$$

From uniform convergence and \eqref{facil} 
we have \begin{align}\label{eq1}
\psi(s, \nu) = \lim_{n \to \infty} \psi(\text{exp}(u_n), \nu_n) \geq \e
\end{align}
Replacing $\nu$ with an ergodic component of $\nu$ satisfying (\ref{eq1}) we can suppose $\nu$ is $s$-ergodic.

Let $p$ denote  the natural projection of $UF $ onto $M^\alpha$ and  let $\mu' = p_{*}\nu$.  Clearly $\mu'$ is  $s$-invariant and ergodic.  We show  that  $\lambda^F_+(s,\mu')$, the average top fiberwise Lyapunov exponent, is positive.  
 Indeed
for $\nu$-almost every $(x_0,v_0)$ in $UF$, it follows from the pointwise ergodic theorem and the chain rule that
\begin{align*}
\e &\leq \int_{UF} \log \bigg( \frac{\|D_xs(v)\|}{\|v\|} \bigg) \ d\nu(x,v) \\
&=\lim_{N\to \infty} \frac{1}{N}
\sum_{k=0}^{N-1} \log \bigg( \frac{\|D_{s^k\cdot x}s(UD_{x_0}s^kv_0)\|}{\|UD_{x_0}s^k v_0\|} \bigg)\\
&=\lim_{N\to \infty} \frac{1}{N}\log \left(  {\|D_{x_0}s^{N}(v_0)\|}.
 \right)
\end{align*}
As
$\inf_{N\to \infty} \frac{1}{N}\log \left(  \restrict{\|D_{x_0}s^{N}}{F} \| \right)\ge \epsilon $
 for $\mu'$-\ae $x_0$, it follows  $\lambda^F_+(s,\mu') \ge \epsilon.$

Finally, averaging $\mu'$ against a \Folner sequence in $A$ and passing to a subsequential limit $\mu$, from Lemma \ref{lemma:averagingiscool}
we have that $\mu$ is $A$-invariant and $\lambda^F_+(s,\mu)\ge  \lambda^F_+(s,\mu')  > 0$.
This completes the proof of Proposition \ref{prop:lyapons}.

\section{Proof of Proposition \ref{prop:Ginvariantmeasure}}
\label{sec:AWBpart}
To prove   Proposition \ref{prop:Ginvariantmeasure}  we apply  an averaging argument to improve certain invariance properties  of the  $A$-invariant measure  on $M^{\alpha}$ with positive exponents produced in Proposition \ref{prop:lyapons}.
Using  measure rigidity results from  homogenous  dynamics, the projection of the averaged measure $\hat \mu$ to $G/\Gamma$ will be the Haar measure.   Using the key technical proposition of \cite{AWBFRHZW-latticemeasure} and the algebraic results in subsection \ref{subsection:algebra}, we deduce that  $\hat \mu$ is in fact $G$-invariant.  We first recall some facts from homogeneous dynamics, particularly a number of results related to Ratner's measure classification theorem,  and  then describe the averaging arguments in the proof.  To illustrate the general argument, the averaging argument is explained for the special case of $\Sl(n,\R)$ in Section \ref{sec:SL}.

\subsection{Facts from homogeneous dynamics}\label{ss:homo}
Let $G$ be a connected, semisimple Lie group  and let $\Gamma\subset G$ be a lattice. Recall that a nilpotent subgroup $U\subset G$ is called \emph{unipotent} if $\ad(u)-\id$ is a nilpotent for every element $u\in U$.
Let $U= \exp \lieu \subset G$ be a unipotent subgroup.
Let $\{b_1, \dots , b_k\}$ be a \emph{regular basis} for $\lieu$ (see \cite{MR1291701}) and for $\vecm= (m_1, \dots, m_k)\in [0,\infty)^k$ let
	$$F_\vecm= \{\exp (t_1 b_1) \cdot \dots \cdot \exp (t_kb_k): 0\le t_j\le m_j \}\subset U.$$
Let $|F_\vecm|$ denote the Haar measure of $F_\vecm$ in $U$.
Recall for $x\in G/\Gamma$ we write  $\nu_x^{F_\vecm} = F_\vecm\ast \delta_x .$  Also recall that a measure $\nu$ is called {\em homogeneous} if there is a closed subgroup $L<G$ such $\nu$ is Haar measure on a closed $L$ orbit in $G/\Gamma$.

\begin{theorem}[Ratner, Shah]\label{thm:ratner}
Let $X= G/\Gamma$ and let $U$ be a unipotent subgroup.    The following hold:
\begin{enumlemma}
	\item \label{ratner1}Every ergodic, $U$-invariant measure is homogeneous \cite[Theorem 1]{MR1262705}.
	\item \label{ratner2} The orbit closure $\orb_x:= \overline{\{u \cdot x :u\in U\} }$ is homogenous for every $x\in G/\Gamma$  \cite[Theorem 3]{MR1262705}.
	\item \label{ratner4} The orbit  ${F_\vecm \cdot x }$ equidistributes in $\orb_x$; that is $\nu_x^{F_\vecm}$ converges to the Haar  measure on $\orb_x$ as $m_1\to \infty, \dots, m_k\to \infty$  \cite[Corollary 1.3]{MR1291701}.
\item \label{ratner3} \label{thisone} Let $A=\exp \liea$ be a maximal split Cartan subgroup, 
	 let $\beta$ be a restricted root of $\lieg$ relative to $\liea$, and  
	 let $\mu$ be a   $G^{[\beta]}$-invariant Borel probability measure on X.  
	If $\mu$ is $A$-invariant then $\mu$ is $G^{[-\beta]}$-invariant.
\end{enumlemma}
\end{theorem}
Note that \ref{thisone} follows from  \cite[Theorem 9]{MR1262705} and the structure of $\mathfrak{sl}(2,\R)$-triples.

Given $x\in G/\Gamma$, let $m^U_x$ denote   the Haar measure on the homogeneous manifold $\orb_x$ in Theorem \ref{thm:ratner}\ref{ratner2}.
Given a measure $\mu$ on $G/\Gamma$ let $$U\ast \mu = \int m^U_x \ d\mu(x).$$
\begin{proposition}\label{thm:averaginghomo}
Let $A=\exp \liea$ be a maximal split Cartan subgroup and let $U= \exp \lieu$ be a unipotent subgroup normalized by $A$.   Let $\mu$ be a Borel probability measure on $G/\Gamma$.  Then
\begin{enumlemma}
	\item \label{averaginghomo1} $ {F_\vecm } \ast \mu\to U\ast \mu$ for any $m_1\to \infty, \dots m_k\to \infty$;
	\item \label{averaginghomo2}If $\mu$ is $A$-invariant, then $U\ast \mu$ is $AU$-invariant;
	\item \label{averaginghomo3}If $\mu$ is $A$-invariant and $A$-ergodic, then  $U\ast \mu$ is $A$-ergodic.
\end{enumlemma}
\end{proposition}

\begin{proof}
For  $x\in G/\Gamma$ we have that $$\nu_x^{F_\vecm }:= F_\vecm \ast \delta_x   $$   converges to the Haar measure $m_x ^U$ on the orbit closure $\orb_x$ of $U\cdot x$.  By dominated convergence we have
$${F_\vecm } \ast \mu = \int \nu_x^{F_\vecm } \ d \mu(x) \to \int m_x^U  \ d \mu(x)=   U\ast \mu$$
and \ref{averaginghomo1} follows.

For \ref{averaginghomo2}, note that if $s\in A$ and if $\{F_\vecm\}$ is a \Folner sequence as above, then $\{s F_\vecm s\inv\}$ is also a \Folner sequence as above.  From the $s$-invariance of $\mu$  and equidistribution in Theorem \ref{thm:ratner}\ref{ratner4} we have that
\begin{align*}
s_* (U\ast \mu)
	&= s_* \left(\lim  \int \nu_x^{F_\vecm } \ d \mu(x) \right)
	\\&
	=\lim s_* \left(  \int \nu_x^{F_\vecm } \ d \mu(x) \right)\\
	&=\lim   \int \nu_{s\cdot x}^{s F_\vecm s\inv  } \ d \mu(x)
	\\&
	=\lim   \int \nu_{ x}^{s F_\vecm s\inv  } \ d \mu(x) \\
	&=    \int  m^U_x  \ d \mu(x).
\end{align*}

For \ref{averaginghomo3}, first write $\calE^U$ for the ergodic decomposition of $U\ast \mu$ for the action of $U$.  By definition, $\calE^U$ coincides with the $(U\ast \mu)$-measurable hull of the partition of $X$ into $U$-orbits.
Let $\{\mu^{\calE^U}_x\}$ denote a family of conditional measures for this partition.   The \Folner sequence $\{F_\vecm\}$ satisfies a pointwise ergodic theorem as $\vecm\to\infty$.  Since each homogeneous measure $m^U_x$ is $U$-ergodic, it follows  for $\mu$-a.e.\ $x'$ and $\mu^{\calE^U}_{x'}$-a.e.\ $x$ that $$\mu^{\calE^U}_x = m^U_{x'}.$$

Let $\phi$ be a bounded, $A$-invariant Borel function.  Using that $U=\exp \lieu$ is unipotent and normalized by $A$, we may select $s_0\in A$ such that $U$ is contracted by $s_0$; that is, $\lieu  \subset \bigoplus_{\beta(s_0)<0} \lieg^\beta.$ By the pointwise ergodic theorem (for the action of $s_0$), $\phi$   coincides modulo $U\ast \mu$ with a  $U$-invariant function.  This follows from the density of uniformly continuous functions in $L^1(U\ast \mu)$.
In particular, the partition into level sets of $\phi$ is coarser (mod $U\ast \mu$) than $\calE^U$.  It follows for $\mu$-a.e.\ $x'\in X$ and $\mu^{\calE^U}_{x'}$-a.e.\ $x\in X$ that $$\phi(x) = \int \phi \ d\mu^{\calE^U}_x= \int \phi \ d m^U _{x'} .$$
In particular, for $(U\ast \mu)$-a.e.\ $x\in X$ there is $x'$ such that $\phi(x) = \int \phi \ d m^U_{x'} .$

Consider the function $\Phi\colon X\to \R$, $$\Phi(x') = \int \phi \ d m^U_{x'}.$$
We have that $\Phi$ is  $\mu$-measurable and, since $U\ast \mu$ is $A$-invariant, $\Phi$ is also $A$-invariant.   From the $A$-ergodicity of $\mu$, $\Phi$ is constant $\mu$-a.e. which implies $\phi$ is constant $\mu_\infty$-a.e.
\end{proof}

\begin{remark}
In the averaging arguments below, we frequently encounter non-ergodic invariant measures $\mu$ on the fiber bundle $M^\alpha$ that project to measures in $G/\Gamma$ with certain desired properties.  To overcome non-ergodicity in our arguments, one may use either Proposition \ref{thm:averaginghomo}\ref{averaginghomo3} or Lemma \ref{lemmatoskipergodicity}.   In the arguments appearing below, we use Lemma \ref{lemmatoskipergodicity} and we never actually use  Proposition \ref{thm:averaginghomo}\ref{averaginghomo3}.    The approach using Proposition \ref{thm:averaginghomo}\ref{averaginghomo3} appears in other versions of the averaging procedure; see for example \cite{BrownTriestinoEtc}.
\end{remark}

\subsection{Averaging argument for $G= \Sl(n,\R)$}\label{sec:SL}
We explain the first step of the proof of Proposition \ref{prop:Ginvariantmeasure} in the case  $G= \Sl(n,\R)$, $n\ge 3$.  Taking the Cartan involution  $\theta\colon \sl(n,\R)\to \sl(n,\R)$ to be $\theta(X) = -X^{\mathrm t}$ we have
$$A =\{ \diag(e^{t_1}, e^{t_2}, \dots, e^{t_n})\}= \left\{\left(\begin{array}{cccc}  e^{t_1}&   &   &      \\  &  e^{t_2} &   &      \\  &   & \ddots  &    \\  &   &   &    e^{t_n} \end{array}\right)\right\}$$
where $  t_1 + t_2 +\dots  + t_n =0.$ Also, $\liec= \{0\}$,
$C$ is the finite group with $\pm 1$ along the diagonals, $K = \So(n)$ and (relative to the standard base)
$$N= \left\{\left(\begin{array}{ccccc}1  &  \ast  & \ast  & \dots  & \ast  \\  & 1  &   \ast& \dots  & \ast  \\  &   &  \ddots &   &  \vdots  \\  &   &   &  1 &  \ast \\  &   &   &   &1  \end{array}\right)\right\}.$$
For $i\neq j\in \{1,\dots n-1\}$ let $\beta_{i,j} \colon A\to \R$ be the linear functional
$$\beta_{i,j}(\diag(e^{t_1}, e^{t_2}, \dots, e^{t_n})) = t_i- t_j.$$
These are the roots of $\sl(n,\R)$ and the    standard base for $\Sigma(\sl(n,\R))$  is $$\Pi=\{\alpha_1 = \beta_{1,2}, \alpha_2 = \beta_{2,3}, \dots , \alpha_{n-1} = \beta_{n-1,n}\}.$$

To prove Proposition \ref{prop:Ginvariantmeasure} it is enough to find an $A$-invariant measure $\mu'$ on $M^\alpha$ with a non-zero fiberwise Lyapunov exponent projecting to the Haar measure on $G/\Gamma$.  By Proposition \ref{prop:LieFacts2} and Proposition \ref{thm:nonresonantimpliesinvariant} below, such a measure will automatically be $G$-invariant.

By the hypotheses of Proposition \ref{prop:Ginvariantmeasure}, we have an ergodic, $A$-invariant measure $\mu$ with a non-zero fiberwise Lyapunov exponent $\lambda^F_\mu\colon A\to \R$.  Note that  $\mu$ need not project to the Haar measure on $G/\Gamma$.  Our goal will be to average $\mu$ over various subgroups of $G$ in order to obtain a new $A$-invariant measure $\mu'$ projecting to the Haar measure.  The subtlety of the argument is to choose the subgroups so that the fiberwise Lyapunov exponents do not vanish after averaging.

Recall that $\lambda^F_\mu\colon A\to \R$ and each $\beta_{i,j}\colon A\to \R$ are non-zero linear functionals.   Consider the linear span $V$ of $\{\alpha_2, \cdots, \alpha_{n-1}\}$ in $\liea^*$.  It may be that $\lambda^F_\mu\in V$.   However, given a permutation matrix (that is, an element of the Weyl group) $P\in \Sl(n,\R)$, let $$P(\lambda^F_\mu)(s) = \lambda^F_\mu(P\inv  s P).$$  One may  check (as the Weyl group acts irreducibly on $\liea^*$) that  $P(\lambda^F_\mu)\notin V$ for some $P$.  Thus, up to conjugating $G$ by a permutation matrix, without loss of generality we may assume  $\lambda^F_\mu\colon A\to \R$  is not in  the linear span   of $\{\alpha_2, \cdots, \alpha_{n-1}\}$.

Let $U$ be the unipotent subgroup
$$
U= \left\{\left(\begin{array}{ccccc}1  &  0  & 0     & \cdots &0 \\  & 1  &   \ast   & \cdots  &\ast \\
 &     &  \ddots &   &  \vdots \\ &  &   &      1 &  \ast \\  &   &   &     &1  \end{array}\right)\right\}$$
and let $$s_1=  \diag\left(\frac{1}{6^{n-1}}, 6, \cdots  , 6\right)\in A.$$
Note that $s_1$ commutes with every element of $U$ and since  $\lambda^F_\mu$ is not in the linear span of  $\{\alpha_2, \cdots, \alpha_{n-1}\}$, $$\lambda^F_\mu(s_1)\neq 0.$$  Replacing $s_1$ with $s_1\inv$, we may assume  $\lambda^F_\mu(s_1)> 0$.

Take a \Folner sequence along $U$ as in Proposition \ref{thm:averaginghomo}, average the measure $\mu$, and pass to a subsequential limit $\mu_1$.  From Proposition \ref{thm:averaginghomo}, we have that $\mu_1$ projects to an $AU$-invariant measure $\hat \mu_1$ in $G/\Gamma$.  Note however that $\mu_1$ may not be $AU$-invariant.  From Lemma \ref{lemma:averagingiscool} however, $\mu_1$ is $s_1$-invariant, $U$-invariant and $\lambda^F_+(s_1,\mu_1)>0$.  Averaging $\mu_1$ along  a \Folner sequence in $A$ and taking  a subsequential limit $\mu_2$, we have $\mu_2$ is $A$-invariant (and in fact { $(AU)$-invariant),}  $\lambda^F_+(s_1,\mu_2 )>0$.  Moreover, as the projection $\hat \mu_1$ of $\mu_1$ is an $AU$-invariant measure,   $\mu_2 $ and $\mu_1$ project to the same $AU$-invariant measure $\hat \mu_1 = \hat \mu_2$ in $G/\Gamma$.  From Theorem \ref{thm:ratner}\ref{ratner3}, it follows  $\hat \mu_1= \hat \mu_2$ is $G'$-invariant where
$$G'= \left\{\left(\begin{array}{ccccc}\ast  &  0 &  0  & \cdots &0 \\  0 & \ast  &  \ast  & \cdots  &\ast \\
0 &\ast  & \ast  &   \  & \ast   \\
\vdots & \vdots &     & \ddots  &  \vdots \\ 0 & \ast  &    \ast   &  \ast &  \ast 
 \end{array}\right)\right\}.$$

Let $\lambda_{+, s_1, \mu_2}\colon A\to \R$ be the linear functional as in Lemma \ref{lemmatoskipergodicity}.  Consider the two roots$$\alpha_1= \beta_{1,2}\colon A\to \R, \quad \quad \delta= \beta_{1,n }\colon A\to \R$$
 (the simple root $\alpha_1$ and the highest root $\delta$.)
Note that $ \lambda_{+, s_1, \mu_2} $ is proportional to at most one of $\alpha_1$ or $\delta$.

 Assume that $ \lambda_{+, s_1, \mu_2} $ is not proportional to  $\alpha_1$.  Let $$U' =  \left\{\left(\begin{array}{ccccc}
 1  &  \ast  & 0     & \cdots &0 \\  & 1  &  0   &    &0 \\
 & &\ddots  &       &  \vdots \\ &  &   &      1 &  0 \\  &   &   &     &1  \end{array}\right)\right\}$$
and select any $s_2 \in \ker \alpha_1\sm \ker \lambda_{+, s_1, \mu_2}$.

Replacing $s_2$ with $s_2\inv$ if necessary, we have  $\lambda^F_+(s_2,\mu_2)\ge  \lambda_{+, s_1, \mu_2}(s_2)>0$.
Average  $\mu_2$ along the one-parameter subgroup $U'$ and pass to a subsequential limit $\mu_3$.   The measure $\mu_3$ projects to an $AU'$-invariant measure $\hat \mu_3$ in $G/\Gamma$.   Average $\mu_3$ along $A$ and pass again to a subsequential limit  $\mu_4$.  We then have
\begin{enumerate}
	\item $\mu_4$ is $A$-invariant;
	\item $\lambda^F_+(s_2,\mu_4)>0$;
	\item $\mu_4$ projects to an $AU'$-invariant measure $\hat \mu_4=\hat \mu_3$  on $G/\Gamma$.
\end{enumerate}
We note that $U'$ commutes with the subgroup $H\subset G'$,
$$H=\left\{ \left(\begin{array}{ccc  cc}1  &  0 & 0    &     \cdots &0 \\  0 & 1  &  0&       &0 \\
0 &\ast  & \ast   &     & \ast   \\
\vdots &   &     &  \ddots     &  \vdots \\ 0 & \ast
 & \ast   &  \ast   & \ast 
 \end{array}\right)\right\}$$

\noindent whence $\hat \mu_3= \hat \mu_4$ is  also invariant under $H$ and $A$.  From Theorem \ref{thm:ratner}\ref{ratner3}, it follows that the projection   $\hat \mu_4= \hat \mu_3$ in  $G/\Gamma$ is invariant under the groups

 $$\left\{
 \left(\begin{array}{ccccc}
 \ast &  \ast  & 0     & \cdots &0 \\  \ast & \ast  &  0   &    &0 \\
 0&0&1&  &0\\
\vdots & &  &      \ddots &  \vdots \\
 0& 0  &0   & \cdots    &1
  \end{array}\right)\right\}, \quad \quad \quad
\left\{  \left(\begin{array}{ccc cc}1  &  0 & 0  &     \cdots &0 \\  0 & \ast  &     \ast  &    &\ast \\
0 &\ast  & \ast  &       & \ast   \\
\vdots &   &   &  \ddots     &  \vdots \\ 0 & \ast
 & \ast   &  \ast &  \ast 
 \end{array}\right)\right\}.$$
 Since these generate $G$, the projection   $\hat \mu_4$ is the Haar measure on  $G/\Gamma$.   Taking an appropriate $A$-ergodic component $\mu'$ of $\mu_4$ we have
 \begin{enumerate}
 \item $\mu'$ is $A$-invariant and $A$-ergodic;
 \item $\mu'$ projects to the Haar measure on $G/\Gamma$;
\item $\lambda^F_+(s_2,\mu')>0$ whence $\mu'$ has a non-zero fiberwise Lyapunov exponent.
 \end{enumerate}

 Above we assumed $ \lambda_{+, s_1, \mu_2} $ was not proportional to $\alpha_1$.  If   $ \lambda_{+, s_1, \mu_2} $ is  proportional to $\alpha_1$ then it is not proportional to $\delta$ and we may take
 $$ U'= \left\{ \left(\begin{array}{ccccc}
 1  & 0  & 0     & \cdots &\ast \\  & 1  &  0   &    &0 \\
 & &\ddots  &       &  \vdots \\ &  &   &      1 &  0 \\  &   &   &     &1  \end{array}\right)\right\}$$
and select any $s_2 \in \ker \delta\sm \ker \lambda_{+, s_1, \mu_2}$.
We may repeat the above arguments (which are now slightly simpler as  $U'$  and $U$ commute) to obtain $\mu_4$ and $\mu'$ with the same properties as before.

\subsection{Averaging argument on $G/\Gamma$}  We present in this and the next subsection the generalization of the averaging procedure described in Section \ref{sec:SL} for general Lie groups. Here, we describe what happens to the projection of the measure to $G/\Gamma$ as we average the measure on $M^{\alpha}$ over various subgroups of $G$.

Let $\lieg$ be a semisimple Lie algebra.  Let $\lieg'$ be a simple ideal of $\lieg$ with rank $\ell\ge 2$ and let $G'\subset G$ be the corresponding analytic  subgroup.  Let $\Sigma$ be the set of  restricted roots of $\lieg'$ and let $\Pi$ be a choice of base generating a system of positive roots $\Sigma_+$.  Let $\Pi = \{\alpha_1, \alpha_2, \cdots, \alpha_\ell\}$ be enumerated such that
$\alpha_1$ is the left-most element in the corresponding Dynkin diagram as drawn in  Appendix \ref{app:Dynkin}.

\begin{proposition}\label{prop:averaginghomo}
With respect to $\Pi$,  let  $\hat\beta$ be either
\begin{enumlemma}
\item $\hat\beta = \delta$, the highest root, if $\lieg'$ is of type $A_\ell$, $B_\ell$,  $D_\ell$, $E_6,$ or $E_7$;
\item $\hat\beta = \delta'$, the 2nd highest root, if $\lieg'$ is of type $C_\ell$, $(BC)_\ell$, $E_8$, $F_4$, or $G_2$.
\end{enumlemma}
Let $\lieu$ be the Lie subalgebra generated by $\{\lieg^{\alpha_2},\dots, \lieg^{\alpha_\ell}\}$ and let $U = \exp \lieu$.
Let $\lieu'$ denote either the Lie subalgebra $\lieg^{\alpha_1}$ or the Lie subalgebra $\lieg^{\hat\beta}$ and let $U' = \exp \lieu'$.

Let $\Gamma\subset G$ be a lattice and let $\mu$ be an $A$-invariant  measure on $G/\Gamma$.  Then
$$U' \ast \left(U \ast \mu\right)$$
is $G'$-invariant.
\end{proposition}

\begin{remark}
The choice of $\hat \beta $  as the highest root $\delta$ or second highest root $\delta'$ in Proposition \ref{prop:averaginghomo} ensures the following two properties hold:
\begin{enumerate}
\item the root subgroups $U^{\hat \beta}$ and $U^{\alpha_j}$ commute for each $2\le j\le \ell$;
\item there is a string of roots $$\beta_0 = \alpha_1 , \  \beta_2, \  \beta_3,  \ \dots, \ \beta_p = \hat \beta$$
such that $\beta_{k} = \beta_{k-1} + \alpha_{j_i}$ for some  $2\le j_i \le \ell$ for each  $1\le i \le m$.
\end{enumerate}
If $\lieg'$ is of type $C_\ell$, $(BC)_\ell$, $E_8$, $F_4$, or $G_2$ the first property holds for the highest root $\hat \beta = \delta$ but the second property fails as $\alpha_1$ has a coefficient of $2$ in $\delta$.  (See Table \ref{table:bull} in Appendix \ref{app:Dynkin}.) The second property is used below to obtain $G'$-invariance after two steps of averaging by obtaining  invariance under root subgroups which generate $G'$; i.e.\ we first pick a $G''$ and average to obtain a $G''$-invariant measure and then make a careful choice of another group to average over to allow us to obtain a $G'$-invariant measure.

Note also in the case that $\Sigma(\lieg')$ is of type $(BC)_\ell$, neither $\hat\beta = \delta'$ nor $\hat\beta = \alpha_1$ is positively proportional to any other root.  In particular $\lieu'= \lieg^{\hat\beta}$ is, in fact, a Lie subalgebra.

\end{remark}

\begin{proof}[Proof of Proposition \ref{prop:averaginghomo}]
Note that $U \ast \mu$ is $U$-invariant.  Let $\nu $ denote $  U' \ast \left(U \ast \mu\right)$.

Consider first the case that $\lieu' = \lieg^{\hat\beta}$.
From the choice of ${\hat\beta}$,  $\lieg^{\hat\beta}$ commutes with each of $\lieg^{\alpha_j}$ for every $2\le j\le \ell$.  From Lemma \ref{lemma:averagingiscool}\ref{averagingiscool1}, $\nu$ is $U$-invariant.  From Proposition \ref{thm:averaginghomo}\ref{averaginghomo2} the measure  $\nu$ is also $A$-invariant.  It follows  from Theorem \ref{thm:ratner}\ref{ratner3} that   $\nu$ is $\exp(\lieg^{-\alpha_j})$-invariant for $2\le j\le \ell$.  From the choice of ${\hat\beta}$ and examining tables of positive roots, there is a sequence of roots $\alpha_1 =   \beta _0 , \beta_1 ,\dots, \beta_p = \hat\beta$ where $\beta_{k-1} = \beta_{k} + (-\alpha_j)$ for some $2\le j\le \ell$ and every $1\le k\le p$.    It follows that $\nu$ is $\exp(\lieg^{\alpha_1})$-invariant.   It then follows that $\nu$ is $G'$-invariant.

In the case that $\lieu' = \lieg^{\alpha_1}$ we first observe that, as $U \ast \mu$ is $U$-invariant,   $U \ast \mu$ is $\exp(\lieg^{-\alpha_j)} $-invariant for every $2\le j\le \ell$.  Since $\lieg^{\alpha_1} $ commutes with  $ \lieg^{-\alpha_j } $ for every $2\le j\le \ell$ it follows that $\nu$ is  $\exp(\lieg^{-\alpha_j)} $-invariant for every $2\le j\le \ell$.  As $\nu$ is $A$-invariant, it follows that $\nu$ is $U$-invariant and, as above, $\nu$ is $G'$-invariant.
\end{proof}

\subsection{Averaging argument on $M^\alpha$} 
Recall that by Remark \ref{rmk:nocompacts} we may assume that  $G$ is a connected, semisimple Lie group with finite center, no compact factors, and all almost-simple factors of real-rank at least 2. Recall the $G$-action on $X= M^{\alpha}$ preserves   the fiberwise tangent  bundle $F= \ker D \pi$. Let $A=\exp \liea  \subset G$  be our fixed maximal split Cartan subgroup.

We assume as in  Proposition \ref{prop:Ginvariantmeasure} and that there is an   $s\in A$ and an $A$-invariant Borel probability measure $\mu$ on $M^\alpha$ with $\lambda_+^F(s,\mu)>0$.
Let $\lieg= \bigoplus_{k=1}^p \lieg'_k$ be the decomposition of $\lieg$ into ideals.  For each $\lieg'_k$, let $G'_k\subset G$ be the corresponding analytic subgroup.  
To complete the proof of  Proposition \ref{prop:Ginvariantmeasure}, we show the following.
\begin{lemma}\label{lem:kkdkd} For $1\le j\le p$,
if the projection of $\mu$ to $G/\Gamma$ is $G'_k$-invariant for all $1\le k\le {j-1}<p$  then there is
an  $s\in A$ and an $A$-invariant Borel probability measure $ \mu'$ on $M^\alpha$ with $\lambda_+^F(s,  \mu')>0$ such that the projection of $\mu'$ to $G/\Gamma$ is
$G'_k$-invariant for all $1\le k\le j\le p$.
\end{lemma}

\begin{proof}

Fix such $G_{j}'$ with Lie algebra $\lieg_{j}'$ and note that $G_j'$ has rank at least $2$.
Let $U, U'$ be as in Proposition \ref{prop:averaginghomo} where the choice of base $\Pi$ and  ${\hat\beta}$ determining $U$ and $U'$ will be made explicit in the proof of Claim \ref{claim:stupid} below.
 Let $F_{i}$,  $F_{i}'$, $F''_{i}$ be \Folner sequences along the nilpotent subgroups $U, U'$, and $A$, respectively, of the type discussed in Section \ref{ss:homo}.  
With $\mu_0 = \mu$, passing to subsequential limits, we may assume we have the following  sequences  of measures  converging in the weak-$*$ topology on $M^\alpha$:
\begin{enumerate}
\item $F_{i_k} \ast \mu_0\to \mu_1$;
\item $F''_{i'_k} \ast \mu_1 \to \mu_2$;
\item $F'_{i''_k} \ast \mu_2 \to \mu_3$;
\item $F''_{i'''_k} \ast \mu_3 \to \mu_4$.
\end{enumerate}
Note that $\mu_2$ and $\mu_4$ are $A$-invariant.
Let $ \mu'  =\mu_4$.
We have the following claim.
 \begin{claim}\label{claim:stupid}
There is a choice of base $\Pi\subset \Sigma(\lieg)$ and a choice of $\hat \beta$ in  Proposition \ref{prop:averaginghomo} such that for  $U$ and $U'$ as in Proposition \ref{prop:averaginghomo}, \Folner sequences $F_{j}$,  $F_{j}'$, $F''_{j}$ as above, and $\mu'$  as above 
\begin{enumlemma}
\item \label{stupid1} $\mu'$ projects to a measure on  $G/\Gamma$ that is $G'_k$-invariant for all $1\le k\le j$;
\item   \label{stupid2}  $\lambda_+^F(s',\mu')>0$ for some $s'\in A$.
\end{enumlemma}
\end{claim}
Lemma \ref{lem:kkdkd}  follows immediately from the above claim.
\end{proof}
We finish the proof of Lemma \ref{lem:kkdkd}  with the proof of Claim \ref{claim:stupid}.
\begin{proof}[Proof of Claim \ref{claim:stupid}]
For any choice of $\Pi$ and choice of $ \hat \beta$,
let $\hat \mu_i$ denote the image of $\mu_i$ in $G/\Gamma$.  We have that $\hat \mu_0$ is $A$-invariant.  We have that $\hat \mu_1= U\ast \hat \mu_0$ is $AU$-invariant whence $\hat \mu_2 = \hat \mu_1$.  From Proposition \ref{prop:averaginghomo} we have that $\hat \mu_3= U'\ast (U\ast \hat \mu_0)$ is $G'_j$-invariant.  As $U\subset G_j'$ and $U'\subset G_j'$ and as $G_k'$ and $G'_{k'}$ commute for $k\neq k'$, it follows from Lemma \ref{lemma:averagingiscool}\ref{averagingiscool1} that
$\hat \mu_3$ is $G_k'$-invariant for all $1\le k\le j-1$.   Then clearly $\hat \mu_4$ is 
$G_k'$-invariant for all $1\le k\le j$.  Conclusion  \ref{stupid1} follows.

For \ref{stupid2} recall that we assume $\lambda_+^F(s,\mu_0)>0$ for some $s\in A$.  Recall the  linear functional  $\lambda_{+, s, \mu_0}\colon A\to \R$  with $\lambda_{+, s, \mu_0}(s) = \lambda_+^F(s,\mu_0)$.
  Also, recall that   restricted roots  $\beta\colon A\to \R$ are linear functionals on $A$.

  We claim there is a choice base
   $\Pi= \{\alpha_1, \dots, \alpha_\ell\}$ so that  $\lambda_{+, s, \mu_0}$ is not in the linear span of  $\{\alpha_2, \dots, \alpha_\ell\}$.
  Indeed, the Weyl group of $\Sigma(\lieg'_j)$ acts irreducibly on $(\liea\cap\lieg'_j)^*$  and simply transitively on bases  $\Pi$ of $\Sigma(\lieg'_j)$.     Moreover the Weyl group preserves angles and lengths so if $\Pi= \{\alpha_1, \dots, \alpha_\ell\}$ is a base of $\Sigma(\lieg'_j)$ and  $\Pi'= \{\alpha_1', \dots, \alpha_\ell'\} = \{ w(\alpha_1), \dots, w(\alpha_\ell)\}$ is the image of $\Pi$ under an element $w$ in the Weyl group, then the vertices  $\{\alpha_1', \dots, \alpha_\ell'\}$ and $\{\alpha_1, \dots, \alpha_\ell\}$ generate the same Dynkin diagram with the same ordering on the vertices.  For a fixed  $\Pi'= \{\alpha_1', \dots, \alpha_\ell'\} $, there is an element $w$ of the Weyl group such that $w(\lambda_{+, s, \mu_0})$ is not in the linear span of $\{\alpha_2', \dots, \alpha_\ell'\} $. Then, letting $\Pi= \{\alpha_1, \dots, \alpha_\ell\}$  map to $\Pi'$ under $w$, we have that $\lambda_{+, s, \mu_0}$ is not in the  linear span of  $\{\alpha_2, \dots, \alpha_\ell\}$.

 We fix this choice of $\Pi= \{ \alpha_1, \dots, \alpha_\ell\}$ for the remainder.

Let $U$ be as in Proposition \ref{prop:averaginghomo} for the above choice of $\Pi$.
Fix   $s_1\in A\sm\ker\lambda_{+, s, \mu_0}$   such that $\alpha_j(s_1) = 0 $ for all $2\le j\le \ell$.  Replacing $s_1$ with $s_1\inv$ if needed we   have
\begin{enumerate}
	\item $ U$ commutes with $s_1$;
	\item $\lambda_+^F(s_1,\mu_0)\ge \lambda_{+, s, \mu_0}(s_1)>0$.
\end{enumerate}
In then follows from Lemma \ref{lemma:averagingiscool} that
\begin{enumerate}
	\item $\mu_1$ is $s_1$-invariant;
	\item $\lambda_+^F(s_1,\mu_1)\ge\lambda_+^F(s_1,\mu_0)>0$;
	\item $\lambda_+^F(s_1,\mu_2)\ge \lambda_+^F(s_1,\mu_1)>0$.
\end{enumerate}

As $\mu_2$ is an $A$-invariant measure on $M^\alpha$, there is a linear functional
 $\lambda_{+, s_1, \mu_2}\colon A\to \R$ with $ \lambda_{+, s_1, \mu_2}(s_1) =\lambda_+^F(s_1,\mu_2) >0$.
Let $\hat \beta$ be as in Proposition \ref{prop:averaginghomo} (relative to the choice of $\Pi$ above).  Note that $\hat \beta$ and $\alpha_1$ are not proportional.  In particular $\lambda_{+, s_1, \mu_2}$ is proportional to at most one of $\{ \hat \beta,  \alpha_1\}$.    Let $\beta'\in \{ \hat \beta,  \alpha_1\}$ be such that $\beta'\neq c\lambda_{+, s_1, \mu_2}$ for any $c\in \R$ and take $\lieu'$  in Proposition \ref{prop:averaginghomo}  to be $\lieu'= \lieg^{\beta'}$.  Fix $s_2\in A$ with $\beta'(s_2 ) = 0$ and $\lambda_{+, s_1, \mu_2}(s_2) >0$.

From Lemma \ref{lemma:averagingiscool} 
we have that
\begin{enumerate}
	\item $\mu_3$ is $s_2$-invariant;
	\item $\lambda_+^F(s_2,\mu_3)\ge\lambda_+^F(s_2,\mu_2)\ge \lambda_{+, s_1, \mu_2}(s_2)>0$;
	\item $\lambda_+^F(s_2,\mu_4)\ge \lambda_+^F(s_2,\mu_3)>0$.
\end{enumerate}
Taking $s'= s_2$ completes the proof of the claim.  \end{proof}

 \subsection{Proof of  Proposition \ref{prop:Ginvariantmeasure}} \label{sec:mainproof}
From Lemma  \ref{lem:kkdkd} it follows that there exists  an  $s\in A$ and an $A$-invariant Borel probability measure $ \mu'$ on $M^\alpha$ with $\lambda_+^F(s,  \mu')>0$ such that the projection of $\mu'$ to $G/\Gamma$ is
$G$-invariant.  In particular, $\mu'$ projects to the Haar measure on $G/\Gamma$.

Let $C$ denote the   centralizer of $A$ in $K$, and $\mu''= C\ast \mu'$.   Since $C$ commutes with $A$ we have
\begin{enumerate}
	\item $\mu''$ is $(CA)$-invariant;
	 \item $\lambda_+^F(s,\mu'')\ge \lambda_+^F(s,\mu')>0$;
	 \item  $\mu''$ projects to the Haar measure on $G/\Gamma$.
\end{enumerate}
Consider a $(CA)$-ergodic component $\bar \mu$ of $\mu''$.  As the Haar measure on $G/\Gamma$ is $(CA)$-ergodic by Moore's ergodicity Theorem, it follows that any such 
$\bar \mu$ projects to the Haar measure on $G/\Gamma$.  With $s$ as above, we may select $\bar\mu $ so that  $\lambda_+^F(s,\bar \mu)>0$.

\begin{definition}Given an $A$-invariant, $A$-ergodic measure $\mu$ on $M^\alpha$,
let $\calL^F= \{ \lambda_j^F\}$ denote the Lyapunov exponent functionals for  the fiberwise derivative cocycle for the measure $\mu$.
We say a restricted root $\beta \in \Sigma(\lieg)$ is \emph{resonant} with the fiberwise exponents of $\mu$ if there is $\lambda_i^F\in \calL^F$ and $c>0$ with $$\beta = c \lambda^F_i.$$
If no such $\lambda_i^F$ and $c$ we say $\beta$ is \emph{non-resonant}.
\end{definition}
 Note that resonance and non-resonance descend to coarse equivalence classes of restricted roots $[\beta]\in \hat \Sigma(\lieg).$

We recall the following key observation from \cite{AWBFRHZW-latticemeasure}.  
\begin{proposition}[{\cite[Proposition 5.1]{AWBFRHZW-latticemeasure}}] \label{thm:nonresonantimpliesinvariant}
Let $\bar \mu$ be   an $A$-invariant Borel probability measure on $M^\alpha$ projecting to the Haar measure on $G/\Gamma$.
Let $\mu$ be an $A$-invariant, $A$-ergodic component of $\bar \mu$.
 Then, given a coarse restricted root $[\beta]\in \hat \Sigma$   that is non-resonant with the fiberwise Lyapunov exponents of $\mu$, the measure $\mu$ is $G^{[\beta]}$-invariant.
\end{proposition}

 Note that the group $C$ acts ergodically  (in fact transitively) on the set of $A$-ergodic components of $\bar \mu$. Moreover, as $C$ commutes with $A$, the group $C$ preserves the Lyapunov exponents for the $A$-action with respect to distinct $A$-ergodic components of $\bar \mu$.  In particular, the set of roots of $\lieg$ that are non-resonant with the fiberwise exponents is constant for almost every  (in fact every) $A$-ergodic component of $\bar \mu$.  Let $\Sigma_{\mathrm{NR}, \bar \mu}$ denote the a.s.\ constant collection of restricted roots of $\lieg$ that are non-resonant with the fiberwise exponents (of ergodic components of $\bar \mu$.)

Let $\lieh\subset \lieg$ be the Lie subalgebra generated by $$\liec \oplus \liea \oplus \bigoplus _{\beta\in \Sigma_{\mathrm{NR}, \bar \mu} } \lieg^{[\beta]}.$$
As there are at most $\dim(M)$ fiberwise Lyapunov exponents it follows that there are at most $\dim(M)$ resonant coarse restricted roots.  It follows that  $\lieh$ has resonant codimension at most $\dim(M)$.   As we assume $\dim(M)\le r(\lieg)$, it  follows from Proposition \ref{prop:LieFacts2} that $\lieh$ is parabolic.

 Let $H\subset G$ be the analytic   subgroup with Lie algebra $\lieh$.  Proposition \ref{thm:nonresonantimpliesinvariant} guarantees that  $\bar\mu$ is $H$-invariant.  We claim $H = G$.  Indeed if $\dim (M) < r(G)$ then $\lieg = \lieh$ follows immediately from the minimality of $r(G)$.
If $\dim (M) = r(g)$ and $H\neq G$ then, as $\lieh$ is  parabolic, we have $$\lieh= \liec \oplus \liea \oplus \bigoplus _{\beta\in \Sigma_{\mathrm{NR}, \bar \mu} } \lieg^{[\beta]}.$$
It follows that every fiberwise Lyapunov exponent is positively proportional with some restricted root $\beta$ with $\lieg^{[\beta]} \cap \lieh = 0$.  In particular, there is an $s\in A$ such that $\lambda^F_i(s)<0$ for every fiberwise Lyapunov exponent $\lambda^F_i\in \calL^F$.
However, in case that the $G$-action preserves a smooth volume in the fibers, the sum of all fiberwise exponents is zero, contradicting the existence of such an $s$.   It thus follows under the hypotheses of Proposition \ref{prop:Ginvariantmeasure}(2) that $\bar \mu$ is $G$-invariant.  This completes the proof of  the proposition.

\begin{remark}
If one wants to obtain a weaker bound in Theorem \ref{theorem:main}, one can replace the argument above using \cite[Proposition 5.1]{AWBFRHZW-latticemeasure} with an easier argument using work of Ledrappier and Young \cite{MR819556}.  This was discovered while this paper was under review and is explained in \cite[Section 8.3]{Cantat}, see also
\cite[Propostions 3 and 4]{BDZ}.  This approach gives similar looking dimension bounds on actions, but with the real rank of $G$ replacing $r(G)$ which is worse in all cases except $\Sl(n, \mathbb{R})$.
\end{remark}

\section{Finding Smooth Metrics}

In this section we prove Theorem \ref{thm:strongTplussubexponential}.  In particular, we  establish the existence of an invariant  Riemannian metric from   uniform subexponential growth of derivatives in conjunction with the strong property (T) of Lafforgue.

\subsection{Lafforgue's strong property (T)}
\label{subsection:strongT}
We recall basic facts about strong property (T).  The reader only interested
in the case of $C^{\infty}$ actions may consider only representations into Hilbert spaces
and ignore the class of Banach spaces $\mathcal{E}_{10}$ introduced in \cite{MR3407190}.  This in fact suffices to prove  theorems for actions by $C^k$ diffeomorphisms on a manifold $M$ when $k=\frac{\dim(M)}{2}+2$.

\begin{definition}
Let $\Gamma$ be a group with a length function $\wl$, $X$ a Banach space and $\pi\colon \Gamma \rightarrow B(X)$.
Given $\epsilon>0$, we say $\pi$ has \emph{$\epsilon$-subexponential norm growth} if there exists a constant $L$ such that
$\|\pi(\gamma)\| \leq Le^{\epsilon\wl(\gamma)}$ for all $\gamma \in \Gamma$.
We say $\pi$ has subexponential norm growth if it has $\epsilon$-subexponential norm growth for all
$\epsilon>0$.
\end{definition}

Given a group $\Gamma$ and a generating set $S$, let $\wl$ be the word length on $\Gamma$.
Here we say a group $\Gamma$ has strong property (T) if it has strong property (T) on Banach spaces for Banach spaces of class $\mathcal E _{10}$ in the quantitative sense of \cite[Section 6]{MR3407190}. In what follows $X$ will denote a Banach space and $B(X)$ will denote the bounded operators on $X$. We will always be considering the operator norm topology on $B(X)$ and we will always mean the operator norm when we write $\|T\|$ for $T \in B(X)$.

%

\begin{definition}
A group $\Gamma$ has \emph{strong property (T)} if there exist a  sequence  of probability measures $\mu_n$ in $\Gamma$ supported in the balls $B(n) = \{\gamma \in \Gamma \mid  l(\gamma) \leq n\}$  such that for every Banach space $X \in \mathcal E _{10}$, there exists a constant $t>0$ such that : For any representation $\pi\colon  \Gamma \rightarrow B(X)$ with $t$-subexponential norm growth
the operators $\pi(\mu_n)$ converge exponentially quickly to a projection onto the space of invariant vectors.  That is,  there exists  $0<\lambda<1$ (independent of $\pi$),  a projection $P \in B(X)$ onto the space of $\Gamma$-invariant vectors, and an $n_0\in \N$ such that $\|\pi(\mu_n)-P\| <  \lambda^n$ for all $n\ge n_0$.
\end{definition}

 We recall the following results obtained from combining results in \cite{MR2423763, MR3407190}:

\begin{theorem}
\label{thm:strongT}
Let $G$ be a connected semisimple Lie group with all simple factors of higher-rank and $\Gamma<G$ a cocompact lattice.
Then $G$ and  $\Gamma$ have strong property (T).
\end{theorem}

\begin{proof}
For the connected Lie group, this is proven explicitly in \cite[Section 6]{MR3407190}.  For the cocompact lattices, this follows from that fact using the proof of \cite[Proposition 4.3]{MR2423763}.  In particular the $\mu_n$ for $\Gamma$ are constructed there explicitly from $\mu_n'$ for $G$ and the properties we desire all follow immediately from this definition since the function $f$ is chosen in $C_C(G)$. A priori, this produces a sequence of measures $\mu_n$ with support in $B(Dn)$ for some fixed number $D$, but by reindexing one can take measures $\mu_n$ supported in $B(n)$. This is not particularly relevant to applications.
\end{proof}

We summarize here some history of strong property $(T)$ and some drift in the definitions of strong property $(T)$.  Lafforgue's original definition only concluded the existence of a self-adjoint projection onto the invariant vectors \cite{MR2423763}. In that paper, Lafforgue introduced strong property $(T)$ and proved that the groups $\Sl(3,F)$ for $F$ any local field, have strong property $(T)$ for representations on Hilbert spaces.  He also noted that this implied strong property $(T)$ on Hilbert spaces for any Lie group containing $\Sl(3, \mathbb R)$ and for cocompact lattice in all such groups.  In subsequent papers, de la Salle and de Laat modified the definition to explicitly include that the projection was a limit of averaging operators defined by measures, but did not assume that the convergence to the limit was exponential \cite{MR3474958,MR3407190}.  In \cite{MR3474958}, de la Salle proved strong property $(T)$ for a much wider class of Banach spaces for $\Sl(3,\mathbb R)$ and in \cite{MR3407190} de Laat and de la Salle proved strong property $(T)$ for both $\Sl(3, \mathbb R)$ and $\Sp(4,\mathbb R)$ and its universal cover for an even wider class of Banach spaces.  These results combined with existing arguments imply strong property $(T)$ for all higher rank simple Lie groups and for their cocompact lattices.  More recently de la Salle has shown that the definition in \cite{MR2423763} and the definition in \cite{MR3474958,MR3407190} are equivalent if one does not necessarily assume that the measures in question are positive \cite{delaSalle}.  It does, however, follows from the proof of  \cite[Theorem 3.9]{delaSalle} that if one has positive measures converging to the projection then there are positive measures converging exponentially quickly to the projection, namely the convolution powers of any measure close enough to the projection.  All existing proofs of strong property $(T)$ explicitly construct sequences of positive measures converging exponentially to a projection \cite{MR2423763,MR3474958,MR3407190}.  While it is not explicitly relevant here, we remark that this is also true of the proof by Liao of strong Banach property $(T)$ for higher rank simple algebraic groups over totally disconnected local fields \cite{MR3190138}.  We also remark that while many of these results extend the class to of Banach spaces satisfying strong property $(T)$ to include some quite exotic Banach spaces, for our purposes it is enough to know the property holds for $\theta$-Hilbertian spaces.

\subsection{Sobolev spaces of inner products}
\label{subsection:sobolev}

To prove Theorem \ref{thm:strongTplussubexponential} from Theorem \ref{thm:strongT}, we need to realize various spaces of $k$-jets of metrics on $M$ as Banach spaces acted on by $\Gamma$.  What follows is a special case of the discussion in \cite[Section 4]{MR2198325} and we refer the reader there for more details and justifications.  Any result stated in this subsection without a reference can be found there.

We will consider the bundle of symmetric two forms on $M$ written as $S^2(TM^*) \rightarrow M$.  The $k$-jets of sections of $S^2(TM^*)$ are
$$ J^k(S^2(TM^*)) \cong \bigoplus_{i=0}^{k} S^i(TM^*){\otimes}S^2(TM^*).$$
A background Riemannian metric on $M$ defines Riemannian metrics on all associated tensor bundles and hence on $J^k(S^2(TM^*))$. There is a natural inclusion
$$C^k(M, S^2(TM^*)) \subset C^0(M, J^k(S^2(TM^*)))$$
 as a closed subspace, but we note that not every section of $J^k(S^2(TM^*)) \rightarrow M$ is  the $k$-jet of a section of $S^2(TM^*)$.  Given a fixed volume form $\omega$, we denote by  $L^p(M, \omega, J^k(S^2(TM)))$   the space of $L^p$ sections of this bundle equipped with  norm defined by
$$\|\sigma\|^p_p= \int_M \|\sigma(m)\|^p d\omega(m).$$
Here the norm inside the integral is defined by the inner product on $S^2(TM^*)_m$ induced by a fixed background Riemannian metric $g$ on $M$.
 Note that, as $M$ is compact, changing the smooth volume $\omega$ or Riemannian metric $g$ gives an isomorphic $L^p$ space and the identity map between any pair of such spaces is bounded. The set of smooth sections of $S^2(TM^*) \rightarrow M$  are naturally included in $L^p(M, \omega, J^k(S^2(TM^*))$.   Let $W^{p,k}(M, \omega, S^2(TM^*))$ be the completion of the set of smooth sections with respect to this norm which we denote $\|\cdot\|_{p,k}$.  Thus
 $$W^{p,k}(M, \omega, S^2(TM^*)) \subset L^p(M, \omega, J^k(S^2(TM^*)))$$
   is a closed subspace.

The following lemma verifies that all the Sobolev spaces  discussed above are in the class $\mathcal{E}_{10}$.
The reader only interested in $C^{\infty}$ actions should consider the case $p=2$ in which all
spaces  discussed above are   Hilbert.

\begin{lemma}
\label{lemma:sobolovspacesinE10}
The Sobolev spaces $W^{p,k}(M, \omega, S^2(TM^*))$ are in the class $\mathcal E_{10}$.
\end{lemma}

\begin{proof}
We use only three facts about $\mathcal E_{10}$: that it contains Hilbert spaces, that the complex interpolation of a space in $\mathcal E_{10}$ with any other space is in $\mathcal E_{10}$, and that $\mathcal E_{10}$ is closed under taking subspaces. This is equivalent to saying that  $\mathcal E_{10}$ contains all $\theta$-Hilbertian spaces. Given any complex Banach space $V$, the spaces $L^p(M,\omega, V)$ is an interpolation spaces of  $L^2(M, \omega, V)$ with $L^{p'}(M,\omega,V )$ for any $p'>p$ and therefore in $E_{10}$. Taking the complexification of $J^k(S^2(TM^*))$ and then passing back to the closed subspace of real valued sections, we  see that $L^p(M,\omega, J^k(S^2(TM^*)))$ is in $\mathcal E_{10}$.  As the class $\mathcal E_{10}$ is closed under taking closed subspaces,  $W^{p,k}(M, \omega, S^2(TM^*))$ is also in $E_{10}$.
\end{proof}

Denote by $C^{k}(M, S^2(TM^*))$ the space of $C^k$ sections of $S^2(TM^*)$.  In the case that $k$ is not integral,
with $l = \lfloor k \rfloor$ and $\lambda = k-l$ elements of
$C^{k}(M, S^2(TM^*)) = C^{l,\lambda}(M, S^2(TM^*))$
are sections of $S^2(TM^*)$ which
are $l$-times differentiable and whose  order-$l$ derivatives are $\lambda$-H\"{o}lder.
We will need the following special case of the Sobolev embedding theorems.

\begin{theorem}
\label{thm:Sobolev}
There is a bounded inclusion $W^{p,l}(M, \omega ,S^2(TM^*)) \subset C^s(M, S^2(TM^*))$ where $s = l - \frac{n}{p}$.
\end{theorem}

As explained in \cite[Section 4]{MR2198325}, this is an easly consequence of the corresponding embedding theorem for domains in $\mathbb R^n$
and the existence of partitions of unity.  We remark that the spaces $W^{p,l}(M, \omega ,S^2(TM^*))$ are defined relative to a   fixed volume form and  metric.  The background volume form and metric need not be preserved.  In our arguments below, the fact that the volume for and metric are not preserved is controlled by the  uniform subexponential growth of derivatives.


\subsection{Proof of Theorem \ref{thm:strongTplussubexponential}}
\label{subsection:findingmetrics}
To construct a $\Gamma$-invariant metric, we first check that the induced action of $\Gamma$ on appropriate Sobolev spaces has subexponential norm growth. Note that $C^k$ actions preserve the class of $C^{k-1}$ Riemmanian metrics, since metrics are defined on the tangent bundle.

\begin{lemma}
\label{lemma:representationhascontrol}
Let $\alpha\colon  \Gamma \rightarrow \Diff^k(M)$ be an action with uniform subexponential growth of derivatives.  Then
the induced representation on $W^{p,k-1}(M,S^2(TM))$ has uniform subexponential norm growth.
\end{lemma}

To prove Lemma \ref{lemma:representationhascontrol}, the key is to see that subexponential growth of the first derivative implies
subexponential growth of all derivatives.  While this is already observed in  \cite{Hurtado}, we include a proof for
completeness.  We recall a special case of \cite[Lemma 6.4]{MR2198325}.  Here given a diffeomorphism  $\phi\colon M\to M$, we write $\|\phi\|_k$ for the norm
of $\phi$ as an operator on $C^k$ vector fields or equivalently $\|\phi\|_k= \sup_{x{\in}M} \|J^k\phi(x)\|$ where $J^k\phi$ is the $k$-jet of $\phi$ or the induced map on $ J^k(TM) \cong \oplus_{i=0}^{k} S^i(TM^*).$

\begin{lemma}\cite[Lemma 6.4]{MR2198325}
\label{lemma:estimateoncompostion} Let
$\phi_1,{\ldots},\phi_n{\in}\Diff^{k}(M)$. Let
$N_k=\max_{1{\leq}i{\leq}n}\|\phi_i\|_k$ and
$N_1=\max_{1{\leq}i{\leq}n}\|\phi_i\|_1$. Then there exists a
polynomial $Q$ depending only on the dimension of $M$ and $k$ such that:
$$\|\phi_1{\circ}{\cdots}{\circ}\phi_n\|_k{\leq}N_1^{kn}Q(nN_k)$$
for every $n{\in}\N$.
\end{lemma}

From this we deduce the following corollary on subexponential growth of higher derivatives.

\begin{corollary}
\label{corollary:highersubexp}
If $\Gamma$ is a finitely generated group, $M$ is a compact manifold and $\alpha\colon\Gamma{\rightarrow}\Diff^k(M)$ has subexponential growth of derivatives then $\alpha$ also has subexponential growth of higher derivatives.  More precisely, subexponential growth of derivatives for $\alpha$ implies that for all $\epsilon>0$ there exists $L_{\epsilon,k}$ such that
$$\|\alpha(\gamma)\|_k \leq L_{\epsilon,k} e^{\epsilon \wl(\gamma)}$$
\noindent for all $\gamma \in \Gamma$.
\end{corollary}

\begin{proof}[Proof of Corolary \ref{corollary:highersubexp}]
We first remark that exponential growth of derivatives is clearly equivalent to the fact that for all $\epsilon >0$ there exists an $n_0$ such that $\|\alpha(\gamma)\|_1 \leq e^{\epsilon \wl(\gamma)}$ for all $\gamma$ with $\wl(\gamma) \geq n_0$.  Applying Lemma \ref{lemma:estimateoncompostion} to words in $\Gamma$ of length $ln_0$ for $l \in \N$, we see that we have for such words
that $\|\alpha(\gamma)\|_k \leq L e^{(k+1)\epsilon \wl(\gamma)}$ where the $L$ and the $k+1$ instead of $k$ are to absorb the polynomial growth into the exponential.  Letting $L' = \sup _{\wl(\gamma) < n_0}\|\alpha(\gamma)\|_k$, by writing all words as products of words of length $kn_0$ and words of length less than $n_0$, we see that $\|\alpha(\gamma)\|_k \leq LL' e^{(k+1)\epsilon \wl(\gamma)}$ for all $\gamma \in \Gamma$.
\end{proof}

\begin{proof}[Proof of Lemma \ref{lemma:representationhascontrol}]

From Corollary \ref{corollary:highersubexp}, we have that for every $\epsilon$ there is an $L$ such that $\|\alpha(\gamma)\|_k < L e^{\epsilon\wl(\gamma)}$.  Up to relabelling $\epsilon$ and $L$ to account for the action on $S^2(TM^*)$, this implies that for  $\sigma \in J^k(M, S^2(TM^*))$, we have a pointwise bound $\|(\alpha(\gamma)_*{\sigma})(x)\| < \|\sigma(\alpha(\gamma)\inv x)\| L e^{\epsilon\wl(\gamma)}$.
This yields the integral bound
$$\int_M\|(\alpha(\gamma)_*{\sigma})(x)\|^p \ d\omega(x) \leq L^p e^{p \epsilon l(\gamma)} \int_M \|\sigma(\alpha(\gamma){\inv}(x)\|^p \ d\omega(x).$$

Write $\Lambda \alpha(\gamma)$ for the Jacobian of derivative of $\alpha(\gamma)$.
Uniform subexponential growth of derivatives   implies that for every $\epsilon>0$ there is an $F>1$ such that
$\frac{1}{F} e^{-n\epsilon l(\gamma)}  \leq \Lambda\alpha(\gamma)(x) \leq F e^{n\epsilon l(\gamma)}$ for every $x \in M$ where $n=\dim(M)$. By change of variable,
$$\int_M \|\sigma(\alpha(\gamma){\inv}(x)\|^p \ d\omega(x) = \int_M \|\sigma(x)\|^p \Lambda\alpha(\gamma)(x) \ d\omega(x)$$
so we have
$$\int_M\|(\alpha(\gamma)_*{\sigma})(x)\|^p \ d\omega(x) \leq FL e^{(p+n)\epsilon l(\gamma)} \|\sigma\|^p_{p,k}.$$
As $\epsilon>0$ was arbitrary, this completes the proof.  
\end{proof}

\begin{proof}[Proof of Theorem \ref{thm:strongTplussubexponential}]
Fix an initial smooth metric $g$.  From  Theorem \ref{thm:strongT} and Lemma \ref{lemma:representationhascontrol},  there exist measures $\mu_n$ supported on $B(n)$ in $\Gamma$ such that $g_n= \pi(\mu_n)g$ converges to an invariant, possibly degenerate, metric $g_{\mathrm{fin}}\in W^{p,k-1}(M,S^2(TM^{*}))$.    In other words $g_{\mathrm{fin}}$ is a non-negative, symmetric, $2$ at each point but  Theorem \ref{thm:strongT}  does not rule out that $g_{\mathrm{fin}}$ is zero on some vector at some point.
  Note  that each $g_n$ is a linear averages of $g$ under the measure  $\mu_n$ on $\Gamma$ and in particular  does not depend on  $p$ or $k$.  Further note that $\| g_n - g_{\mathrm{fin}} \|_{p,k} \leq C_{p,k}^n$ for some $O < C_{p,k} < 1$ and all $n$ sufficiently large.  Applying Theorem \ref{thm:Sobolev}, it follows that $g_{\mathrm{fin}}$ is in $C^{k-1-\frac{\dim(M)}{p}}$ for all choices of $p$ and is thus $C^{k-1-\beta}$ for all $\beta>0$. If the   action is by $C^{\infty}$ diffeomorphisms, this proves $g_{\mathrm{fin}}$ is $C^{\infty}$.  If the action is $C^2$, the metric  $g_{\mathrm{fin}}$  is only H\"{o}lder.

It remains to check that  $g_{\mathrm{fin}}$   is not degenerate.  This follows
as the averaged metrics $g_n$ degenerate subexponentially while the convergence $g_{\mathrm{fin}}$   is exponentially fast.   To see this explicitly, we check that  $g_{\mathrm{fin}}(v,v)>0$ for any unit vector $v$ in $TM_m$.  The Sobolev embedding theorems imply that $\| g_n - g_{\mathrm{fin}} \|_{0} < K C^n$ for some $0<C<1$, $K>0$, and all sufficiently large $n$.  Choose  $\epsilon>0$ with $Ce^{\epsilon}<1$.  Uniform subexponential growth of derivatives implies that there is a constant $L>0$ such that
$$\|g(D\alpha(\gamma)(v),D\alpha(\gamma)(v))\| \geq L e^{-\epsilon l(\gamma)}.$$
 This implies that
$$g_n(v,v) \geq Le^{-\epsilon n}\|v\|^2.$$
If $g_{\mathrm{fin}}(v,v)=0$
then it would follow that  $g_n(v,v) \leq C^n$ whence $L e^{-\epsilon n} < K C^n$ for all sufficiently large $n$.  But then
$$\frac{L}{K}\leq (Ce^{\epsilon})^n$$
 for all sufficiently large $n$, a contradiction.
\end{proof}

\newpage
\appendix
\section{Table of root data}\label{app:Dynkin}
 The following table includes Dynkin diagrams of all irreducible root systems and an enumeration of the simple roots relative to a choice of  base $\Pi$.  We also include the highest and second highest roots $\delta$ and $\delta'$ relative to the base $\Pi$ and the resonant codimension of all maximal parabolic subalgebras $\lieq_j:= \lieq_{\Pi\sm \{\alpha_j\}}.$
 \DynkTable
\bibliographystyle{AWBmath}
\bibliography{bibliography}
\end{document}